\documentclass[a4paper,12pt]{article}

\usepackage{amsmath,amssymb,amsthm}

\newcommand{\BC}{\mathbb{C}}
\newcommand{\BZ}{\mathbb{Z}}
\newcommand{\BB}{\mathbb{B}}

\newcommand{\CB}{\mathcal{B}}
\newcommand{\CL}{\mathcal{L}}
\newcommand{\CS}{\mathcal{S}}
\newcommand{\CJ}{\mathcal{J}}
\newcommand{\CH}{\mathcal{H}}

\newcommand{\Fh}{\mathfrak{h}}
\newcommand{\FS}{\mathfrak{S}}

\newcommand{\bc}{{\bf c}}
\newcommand{\bx}{{\bf x}}
\newcommand{\by}{{\bf y}}

\newcommand{\bV}{{\bf V}}
\newcommand{\bW}{{\bf W}}
\newcommand{\vm}{{\bf 1}}

\newcommand{\ds}{\displaystyle}
\newcommand{\stm}{\setminus}

\newcommand{\ha}[1]{\widehat{#1}}

\newcommand{\oi}{\omega_{r}^{ii}}
\newcommand{\oij}{\omega_{r}^{ij}}
\newcommand{\li}{L_{r}^{ii}}
\newcommand{\lj}{L_{r}^{jj}}
\newcommand{\lij}{L_{r}^{ij}}
\newcommand{\lab}[2]{L_{r}^{#1}(#2)}

\newcommand{\Lr}{\CL_{r}}
\newcommand{\Vam}{V_{\CJ}}
\newcommand{\Mam}{M_{r}}
\newcommand{\Lon}{\CL_{r}^{(1)}}
\newcommand{\Mon}{M_{r}^{(1)}}

\newcommand{\sm}{\sigma}

\newcommand{\vsp}{\vspace{1.5mm}}

\DeclareMathOperator{\End}{End}
\DeclareMathOperator{\Id}{id}
\DeclareMathOperator{\Sym}{Sym}
\DeclareMathOperator{\Hom}{Hom}
\DeclareMathOperator{\ad}{ad}
\DeclareMathOperator{\Span}{Span}
\DeclareMathOperator{\sn}{sign}
\DeclareMathOperator{\sgn}{sgn}

\theoremstyle{plain}
\newtheorem{thm}{Theorem}[section]
\newtheorem{lem}[thm]{Lemma}
\newtheorem{prn}[thm]{Proposition}
\newtheorem{cor}[thm]{Corollary}

\newtheorem{claim}{Claim}[thm]

\newtheorem*{claim*}{Claim}
\newtheorem*{thm*}{Theorem}

\theoremstyle{definition}
\newtheorem{dfn}[thm]{Definition}

\theoremstyle{remark}
\newtheorem{rmk}[thm]{Remark}
\newtheorem{eg}[thm]{Example}

\makeatletter
\@addtoreset{equation}{section}

\renewcommand\section{\@startsection{section}{1}{0pt}
{-3.5ex plus -1ex minus -.2ex}{1.0ex plus .2ex}{\large\bf}}
\renewcommand\subsection{\@startsection{subsection}{1}{0pt}
{2.5ex plus 1ex minus .2ex}{-1em}{\bf}}
\makeatother

\begin{document}

\setlength{\baselineskip}{18pt}

\title{{\Large \bf
Simplicity of a vertex operator algebra
whose Griess algebra is 
the Jordan algebra of symmetric matrices}}

\author{Hidekazu Niibori\\ 
\small Graduate School of Pure and Applied Sciences, University of Tsukuba,\\
\small Tsukuba, Ibaraki 305-8571, Japan\
(e-mail: {\tt niibori@math.tsukuba.ac.jp})\\[2mm]
and\\[2mm]
Daisuke Sagaki\\ 
\small Institute of Mathematics, University of Tsukuba,\\
\small Tsukuba, Ibaraki 305-8571, Japan\
(e-mail: {\tt sagaki@math.tsukuba.ac.jp})}
\date{}
\maketitle
%
%
\begin{abstract}
Let $r \in \BC$ be a complex number, and 
$d \in \BZ_{\ge 2}$ a positive integer greater than or equal to $2$.
Ashihara and Miyamoto introduced a vertex operator algebra $\Vam$ of 
central charge $dr$, whose Griess algebra is isomorphic to 
the simple Jordan algebra of symmetric matrices of size $d$. 
In this paper, we prove that the vertex operator algebra 
$\Vam$ is simple if and only if $r$ is not an integer.
Further, in the case that $r$ is an integer 
(i.e., $\Vam$ is not simple), we give a generator system of 
the maximal proper ideal $I_{r}$ of the VOA $\Vam$ explicitly. 
\end{abstract}
%
%
\section{Introduction.}
\label{sec:intro}

Let $V=\bigoplus_{n \in \BZ_{\ge 0}} V_{n}$ 
be a vertex operator algebra (VOA for short) 
over the field $\BC$ of complex numbers 
with $Y(\cdot\,,\,z):V \rightarrow (\End V) [[z,z^{-1}]]$ 
the vertex operator. 
As usual, for each $v \in V$, 
we define $v_{m} \in \End V$, $m \in \BZ$, by: 
$Y(v,\,z)=\sum_{m \in \BZ} v_{m}z^{-m-1}$. 
It follows from the axiom of a VOA that 
$V_{2}$ becomes a $\BC$-algebra with the product given by
$u \cdot v:=u_{1}v \in V_{2}$ for $u,\,v \in V_{2}$. 
In addition, we know (see \cite[\S10.3]{flm} and also \cite[\S5]{m}) that 
if $\dim V_{0}=1$ and $\dim V_{1}=0$, then 
the $\BC$-algebra $V_{2}$ is commutative 
(but not necessarily associative), 
which we call the Griess algebra of $V$. 

Various kinds of commutative $\BC$-algebras appear as 
the Griess algebras of VOAs.
The Griess algebra of the moonshine VOA $V^{\natural}$ is isomorphic to 
the 196884-dimensional, commutative $\BC$-algebra introduced by Griess \cite{g}, 
whose automorphism group is isomorphic to 
the Monster sporadic simple group (see \cite[p.319]{flm}); 
the name ``Griess algebra'' is derived from this fact. 
Also, for a given associative, commutative $\BC$-algebra $A$ 
equipped with an $A$-invariant bilinear form, 
Lam \cite{lam1} constructed a VOA whose Griess algebra is 
isomorphic to $A$. 
Further, in \cite{lam2}, \cite{am}, \cite{ash}, 
they constructed some VOAs whose Griess algebras are (simple) 
Jordan algebras; for the definition of a Jordan algebra, 
see \S\ref{subsec:jordan} below. 
In this paper, we will mainly treat VOAs 
introduced in \cite[\S4.1]{lam2} and \cite{am}, 
whose Griess algebras are isomorphic to 
the simple Jordan algebra $\Sym_{d}(\BC)$ of 
symmetric matrices of size $d \in \BZ_{\ge 2}$ 
with entries in $\BC$. 
Because the Jordan algebra $\Sym_{d}(\BC)$ has a strong
connection to symmetric cones and zeta functional 
equations (see \cite{fk}), Ashihara, Miyamoto (see \cite[Introduction]{am}), 
and the authors of this paper expect that 
the results in \cite{am} and this paper 
contribute a VOA theoretical approach to 
the theory of symmetric cones and zeta functional equations. 

An essential difference between 
the VOA introduced in \cite[\S4.1]{lam2} and 
the one introduced in \cite{am} (which we denote by $\Vam$) 
is their central charges. 
The central charge of the former is equal to 
the (fixed) positive integer $d \in \BZ_{\ge 2}$, 
the size of symmetric matrices. On the other hand, 
the central charge of the later is equal to $dr$, 
where $r$ is an arbitrary complex number 
(see Theorem~\ref{thm:am} below). 
In general, the structure and the representation theory of 
a VOA deeply depends on its central charge. For example, 
it is well-known that the simplicity of a Virasoro VOA 
$M_{c,0}/\langle L_{-1}1 \rangle$ (with notation in \cite{w}) and 
the rationality of a simple Virasoro VOA $V_{c}$ 
(with notation in \cite{w}) depend on 
their central charges $c$ (see \cite{w} and also \cite{dmz}). 
Moreover, rational Virasoro VOAs 
(i.e., $V_{c}$ of special central charge $c$, such as $V_{1/2}$) and 
their irreducible modules play very important roles 
in the theory of VOAs. 
So, as in the case of Virasoro VOAs, 
it is quite natural and important to study how 
the VOA $\Vam$ introduced in \cite{am} depends on 
its central charge $dr$. 
In this paper, we study the condition of 
$r \in \BC$ for the VOA $\Vam$ to be simple. 
The main result of this paper is the following theorem, 
which means that the simplicity of $\Vam$ 
also depends on its central charge. 

\begin{thm*}
Keep the notation above. 
The VOA $\Vam$ is simple 
if and only if $r$ is not an integer. 
\end{thm*}

Many important VOAs are obtained as 
the nontrivial simple quotients of 
nonsimple VOAs (e.g., rational Virasoro VOAs, and 
VOAs associated to integrable highest weight 
modules over affine Lie algebras).  
So we are interested in the simple quotient $\Vam/I_{r}$ 
with $r \in \BZ$ rather than the VOA $\Vam$. 
When we study the structure of $\Vam/I_{r}$ and 
its representation theory, it is very important 
to determine some relations in 
$\Vam/I_{r}$ induced by the maximal proper ideal $I_{r}$. 
In this paper, as a first step for studying 
the simple VOA $\Vam/I_{r}$, 
we will give a generator system of 
the maximal proper ideal $I_{r}$ of 
the VOA $\Vam$ explicitly, whose elements are 
singular vectors for a certain Lie algebra $\Lon$, and 
have a high symmetry (see \S\ref{sec:step4} below).

\vsp

This paper is organized as follows. 
In \S\ref{sec:review}, we recall the definition of 
the Griess algebra of a VOA, the definition of a Jordan algebra, and 
the construction of the VOA $\Vam$ introduced by Ashihara and 
Miyamoto \cite{am}. Then we state our main theorem 
(Theorem~\ref{thm:main}), and the plan how we prove 
the theorem. 
In \S\ref{sec:step1}\,--\,\S\ref{sec:step4}, 
following the plan, we will prove some key propositions 
(Propositions~\ref{V=M}, \ref{prop:step1}, \ref{prop:step2}, 
 \ref{prop:step3}, \ref{prop:step4}); 
our main theorem follows immediately from these propositions. 
In \S\ref{sec:step4}, we also give a generator system of 
the maximal ideal $I_{r}$ of the VOA $\Vam$ explicitly 
when $r$ is an integer, that is, when $\Vam$ is not simple. 

\paragraph{Acknowledgments.}
The authors would like to express their sincere gratitude to 
Professor Masahiko Miyamoto, Professor Toshiyuki Abe, 
Professor Hiroki Shimakura, Professor Hiroshi Yamauchi, and 
Dr. Takahiro Ashihara for valuable discussions.

%
\section{Vertex operator algebra whose Griess algebra is a Jordan algebra.}
\label{sec:review}
%
%
\subsection{Griess algebras.}
\label{subsec:griess}
Let
$\bigl(
 V=\bigoplus_{n \in \BZ_{\ge 0}}V_n,\ 
 Y(\cdot\,,z),\ 
 \vm,\ \omega \bigr)$
be a vertex operator algebra (VOA for short),
with $Y(\cdot\,,z):V \rightarrow (\End V) [[z,z^{-1}]]$ 
the vertex operator, $\vm \in V_0$ the vacuum element,
and $\omega \in V_2$ the Virasoro element
(for the details about VOAs, see, e.g., \cite{ll}). 
As usual, for each $v \in V$, 
we define $v_{m} \in \End V$, $m \in \BZ$, by: 
$Y(v,\,z)=\sum_{m \in \BZ} v_{m}z^{-m-1}$. 
For $a, b \in V_2$, we define $a \cdot b:=a_1b$. 
Then it follows from the axiom of a VOA
that $a \cdot b \in V_2$ for every $a, b \in V_2$,
i.e., $V_2$ becomes a $\BC$-algebra with $\cdot$ the product.
In addition, if $V_0=\BC \vm$ and $V_1=\{0\}$,
then the $\BC$-algebra $V_2$ is commutative
(see \cite[\S10.3]{flm} and also \cite[\S5]{m}).
In this case, we call $V_2$ the Griess algebra of $V$.
Note that the Griess algebra of a VOA is 
not necessarily associative.
%
%
\subsection{Jordan algebras.}
\label{subsec:jordan}
Let us recall the definition of a Jordan algebra.
For the details about Jordan algebras, see, e.g., 
\cite{alb46}, \cite{alb48}, and \cite{ja}.
\begin{dfn}
Let $J$ be a $\BC$-algebra 
with the product $a \cdot b \ (a, b \in J)$. 
The $\BC$-algebra $J$ is called a Jordan algebra if 
$a \cdot b=b \cdot a$ and 
$a^2\cdot(b\cdot a)=(a^2\cdot b)\cdot a$ hold
for every $a,b \in J$.
\end{dfn}
Let $\Sym_d(\BC)$ be the set of symmetric matrices of 
size $d \in \BZ_{\ge 2}$ with entries in $\BC$.
It is well-known that 
$\Sym_d(\BC)$ becomes a (simple) Jordan algebra,
where the product is given by: $A \cdot B=
\frac{1}{2}(AB+BA)$ for $A,\,B \in \Sym_d(\BC)$
(see also \cite[Theorem~6\,B]{lam2}).
%
%
\subsection{VOA $\Vam$.}
\label{subsec:am}
Let (and fix) $d \in \BZ_{\ge 2}$. 
In this subsection, we recall a VOA $\Vam$ 
introduced by Ashihara and Miyamoto \cite{am}, 
whose Griess algebra is isomorphic to the Jordan algebra 
$\Sym_d(\BC)$ of symmetric matrices.

Let $\ha{\Fh}$ be an (infinite-dimensional) vector space over $\BC$ 
with a linear basis $\bigl\{v^{i}(m) \mid 1 \le i \le d, \, m \in \BZ \bigr\} \cup \{\bc\}$, 
and define a Lie bracket on $\ha{\Fh}$ by: 
%
%
\begin{equation} \label{hei}
\begin{array}{l}
[v^{i}(m),\ v^{j}(n)] 
  = \delta_{m+n,0}\,\delta_{i,j}\,m \bc 
  \qquad \text{for $1 \le i,\,j \le d$ and $m,n \in \BZ$}, \\[3mm]
[\bc,\ \ha{\Fh}]=\{0\}.
\end{array}
\end{equation}
Denote by $U(\ha{\Fh})$ the universal enveloping algebra of 
the Lie algebra $\ha{\Fh}$, and 
let $U(\ha{\Fh})/\langle \bc-1 \rangle$ be the quotient algebra 
with respect to the two-sided ideal $\langle \bc-1 \rangle$
of $U(\ha{\Fh})$ generated by $\bc-1 \in U(\ha{\Fh})$.
We define a subspace $\CL$ of $U(\ha{\Fh})/\langle \bc-1 \rangle$ as follows:
For $1 \le i,\,j \le d$ and $m,\,n \in \BZ$, we set 
%
%
\begin{equation} \label{eq:dvij}
v^{ij}(m,\,n):= 
 v^{i}(m)v^{j}(n) \mod \langle \bc-1 \rangle. 
\end{equation}
Let
\begin{align*}
& \CB := 
  \left\{v^{ii}(m,\,n) \mid
  \text{$1 \le i \le d$ and 
         $m,\,n \in \BZ$ with $m \le n$} \right\} \ \cup \\
& \hspace{40mm}
  \left\{ v^{ij}(m,\,n) \mid
  \text{$1 \le i < j \le d$ and $m,\,n \in \BZ$} \right\}.
\end{align*}
Then it follows from the Poincar\'e-Birkhoff-Witt theorem that
$\CB \cup \{1 \in U(\ha{\Fh})/\langle \bc-1 \rangle\}$ is 
a linearly independent subset of $U(\ha{\Fh})/\langle \bc-1 \rangle$.
We set
\begin{equation*}
\CL:=\bigl(\Span_{\BC} \CB\bigr) \oplus \BC
\subset U(\ha{\Fh})/\langle \bc-1 \rangle.
\end{equation*}
%
%
\begin{rmk} \label{rem:CB}
Let $1 \le i,\,j \le d$, and $m,\,n \in \BZ$. 
It can be easily seen from the definition \eqref{hei} of
the Lie bracket on $\ha{\Fh}$ that
%
%
\begin{equation} \label{eq:CB}
\begin{array}{l}
v^{ij}(m,\,n)=v^{ji}(n,\,m) \quad
  \text{if $i \ne j$, or if $i=j$ and $m \ne -n$}, \\[3mm]
v^{ii}(m,\,-m)=v^{ii}(-m,\,m)+m.
\end{array}
\end{equation}
In particular, $v^{ij}(m,\,n) \in \CL$
for all $1 \le i,\,j \le d$ and $m,\,n \in \BZ$. 
\end{rmk}

We see by direct computation 
that $[x,y]=xy-yx$ is contained in $\CL$ for every $x,\,y \in \CL$,
and hence $\CL$ becomes a Lie algebra with respect 
to the natural Lie bracket.
Now, let (and fix) $r \in \BC$ be an (arbitrary) complex number.
For each $x,\,y \in \CL$, we define 
%
%
\begin{equation} \label{nLp}
[x,\,y]_{r}:=\pi_{1}([x,y])+r\pi_{2}([x,y]), 
\end{equation}
where $\pi_1:\CL \twoheadrightarrow \Span_{\BC} \CB$
and $\pi_2:\CL \twoheadrightarrow \BC$ denote projections from $\CL$
onto $\Span_{\BC} \CB$ and $\BC$, respectively. 
Then we know from \cite[\S2.1]{am} that 
$[\cdot\,,\,\cdot]_{r}$ is a Lie bracket on $\CL$. 
Let us denote by $\Lr$ the Lie algebra $\CL$ 
with the new Lie bracket $[\cdot\,,\,\cdot]_{r}$. 
%
%
\begin{eg} \label{eg:bracket}
As an example, let us compute
$[v^{ii}(m,\,n),\,v^{ii}(-n,\,-m)]_{r}$ 
for $1 \le i \le d$ and $m,\,n \in \BZ_{> 0}$ 
with $m \le n$. By direct computation and \eqref{eq:CB}, 
we see that in $\CL$, 
\begin{align*}
& [v^{ii}(m,\,n),\,v^{ii}(-n,\,-m)] \\
& \hspace*{10mm}
  = n(1+\delta_{m,n})v^{ii}(-m,\,m)+
    m(1+\delta_{m,n})v^{ii}(-n,\,n)+mn(1+\delta_{m,n}).
\end{align*}
Hence, by the definition \eqref{nLp} of 
the Lie bracket $[\cdot\,,\,\cdot]_{r}$, we obtain 
%
%
\begin{align}
& [v^{ii}(m,\,n),\,v^{ii}(-n,\,-m)]_{r} \notag \\
& \hspace*{5mm} =
  \pi_{1}([v^{ii}(m,\,n),\,v^{ii}(-n,\,-m)])+
  r\pi_{2}([v^{ii}(m,\,n),\,v^{ii}(-n,\,-m)]) \notag \\
& \hspace*{5mm} =
  n(1+\delta_{m,n})v^{ii}(-m,\,m) + 
  m(1+\delta_{m,n})v^{ii}(-n,\,n)+rmn(1+\delta_{m,n}). \label{eq:LB01}
\end{align}
\end{eg} 

We now set
\begin{align*}
\CB_{+} 
  & :=\left\{ v^{ij}(m,\,n) \in \CB \mid 
      m \in \BZ_{\geq 0} \text{ or } n \in \BZ_{\geq 0} \right\}, \\
\CB_{-}
  & := \left\{ v^{ij}(m,\,n) \in \CB \mid m,\,n \in \BZ_{<0}\right\},
\end{align*}
and $\Lr^+:=\bigl(\Span_{\BC} \CB_{+}\bigr) \oplus \BC$. 
It is easily seen that
$\Lr^+$ is a Lie subalgebra of $\Lr$.
Let $\BC \vm$ be a one-dimensional $\Lr^{+}$-module
such that $x \cdot \vm =0$ for all $x \in \CB_{+}$, and 
$s \cdot \vm = s\vm$ for each $s \in \BC \subset \Lr^{+}$. 
Denote by $\Mam$ the $\Lr$-module induced from
the $\Lr^+$-module $\BC \vm$, that is, 
\begin{equation*}
\Mam:=U(\Lr) \otimes_{U(\Lr^+)}\BC\vm.
\end{equation*}

Here we give a linear basis $\BB$ of $\Mam$ by using 
the Poincar\'e-Birkhoff-Witt theorem. 
Take (and fix) a total ordering $\succ$ on the set $\CB_{-}$.
Denote by $\CS$ the set of finite sequences of elements of $\CB_{-}$
that is weakly decreasing with respect to the total ordering $\succ$. 
For $\bx=(x_{p} \succeq x_{p-1} \succeq \cdots \succeq x_1) 
\in \CS$ with $x_{q} \in \CB_{-}$ for $1 \leq q \leq p$,
we set 
\begin{equation*}
w(\bx):=x_{p}x_{p-1} \cdots x_{1}\vm \in \Mam.
\end{equation*}
In view of the Poincar\'e-Birkhoff-Witt theorem, 
$\BB:=\left\{ w(\bx) \mid \bx \in \CS \right\}$ is 
a linear basis of the $\Lr$-module $\Mam$.
%
%
\begin{rmk} \label{rem:pbw}
We see from the definition \eqref{nLp} of 
the Lie bracket on $\CL_{r}$ that $xy=yx$ 
for all $x,\,y \in \CB_{-}$. Therefore, 
if $y_{1},\,y_{2},\,\dots,\,y_{p} \in \CB_{-}$, then 
$y_{1} y_{2} \cdots y_{p} \vm = w(\bx) \in \BB$,
where $\bx \in \CS$ is the sequence of length $p$ 
obtained by arranging $y_{1},\,y_{2},\,\dots,\,y_{p}$ 
in the weakly decreasing order with respect to 
the total ordering $\succ$. 
Also, we note that if $m,\,n \in \BZ_{< 0}$, then 
$v^{ij}(m,\,n)=v^{ji}(n,\,m)$ for every 
$1 \le i,\,j \le d$ (see Remark~\ref{rem:CB}). 
\end{rmk}

If $\bx=(x_{p} \succeq x_{p-1} \succeq \cdots \succeq x_1) \in \CS$ 
with $x_q=v^{i_qj_q}(m_q,\,n_q) \in \CB_{-}$ for $1 \le q \le p$, 
then we define the degree of $w(\bx) \in \BB$ by:
\begin{equation*}
\deg (w(\bx))=-\sum_{q=1}^{p} (m_q+n_q) \in \BZ_{\ge 0}. 
\end{equation*}
Then the $\Lr$-module $\Mam$ admits 
the degree space decomposition as follows: 
\begin{equation*}
\Mam=\bigoplus_{n \in \BZ_{\geq 0}}(\Mam)_n, \quad
\text{where } 
(\Mam)_n:=\Span_{\BC} 
  \bigl\{b \in \BB \mid \deg b = n \bigr\}.
\end{equation*}
Note that 
%
%
\begin{equation} \label{wt01}
(\Mam)_0=\BC\vm, \quad \text{and} \quad 
(\Mam)_1=\{0\}. 
\end{equation}

Define an operator $\lij(m) \in \End(\Mam)$ 
for $1 \leq i,\,j \leq d$ and $m \in \BZ$ by:
%
%
\begin{equation} \label{dLij}
\lij(m)=
 \begin{cases}
 \dfrac{1}{2}\ds{\sum_{h \in \BZ}} v^{ij}(m-h,\,h) 
 & \text{if $i \ne j$ or $m \ne 0$}, \\[8mm]
 \dfrac{1}{2}v^{ii}(0,\,0) + 
 \ds{\sum_{h \in \BZ_{>0}}} v^{ii}(-h,\,h) 
 & \text{if $i=j$ and $m=0$},
 \end{cases}
\end{equation}
and set 
\begin{equation*}
\oij:=\lij(-2)\vm \in (\Mam)_{2},
\quad \text{and} \quad
\omega:=\sum_{i=1}^d \oi \in (\Mam)_{2}.
\end{equation*}
Remark that $\lij(m)=L^{ji}_{r}(m)$ 
for every $1 \le i,\,j \le d$ and $m \in \BZ$ 
(see Remark~\ref{rem:CB}), and hence that 
$\oij=\omega_{r}^{ji}$ for every $1 \le i,\,j \le d$. 
Let $\CJ:=
 \bigl\{\oij,\ \vm \mid 1 \leq i,\,j \leq d \bigr\}
 \subset \Mam$,
and let $\Vam$ be the subspace of 
$\Mam$ spanned by all elements of the form:
$L_r^{i_1j_1}(m_1)L_r^{i_2j_2}(m_2)\cdots L_r^{i_pj_p}(m_p)\vm$ 
with $p \ge 0$, and $1 \leq i_q,\,j_q \leq d$, 
$m_q \in \BZ$ for $1 \leq q \leq p$.
Then, $\Vam$ also admits the degree space decomposition 
induced from that of $\Mam$, i.e., 
$\Vam=\bigoplus_{n \in \BZ_{\geq 0}}(\Vam)_n$
with $(\Vam)_n:=\Vam \cap (\Mam)_n$ 
for $n \in \BZ_{\geq 0}$. We should remark that 
$(\Vam)_0=\BC\vm$ and $(\Vam)_1=\{0\}$ by \eqref{wt01}.

Define a map 
$Y_{0}(\cdot\,,\,z):\CJ \rightarrow \End (\Vam)[[z,z^{-1}]]$ by: 
%
%
\begin{equation} \label{eq:y0}
\begin{array}{l}
Y_{0}(\oij,\,z)=\ds{\sum_{m \in \BZ}} \lij(m) z^{-m-2} 
  \quad \text{for $1 \leq i,\,j \leq d$}, \\[7mm]
Y_{0}(\vm,\,z)=\Id_{\Vam}.
\end{array}
\end{equation}
The following theorem is the main result of \cite{am}.
%
%
\begin{thm} \label{thm:am}
Keep the notation above. 
The map $Y_0(\cdot\,, z):\CJ \rightarrow \End(\Vam)[[z,z^{-1}]]$ 
can be uniquely extended to a linear map
$Y(\cdot\,, z):\Vam \rightarrow \End(\Vam)[[z,z^{-1}]]$ 
in such a way that the quadruple
$\left(
 \Vam=\bigoplus_{n \in \BZ_{\geq 0}}(\Vam)_n,\ 
 Y(\cdot\,, z),\ 
 \vm,\ \omega \right)$ becomes a VOA of central charge $dr$, 
with $\vm$ the vacuum element, and $\omega$ the Virasoro element. 
Furthermore, the Griess algebra of $\Vam$ is 
isomorphic to the Jordan algebra $\Sym_{d}(\BC)$ of 
symmetric matrices. 
\end{thm}

The purpose of this paper is to determine 
the condition of $r \in \BC$ for the VOA $\Vam$ 
to be simple. 
%
%
The following theorem is the main result of this paper. 
%
%
\begin{thm} \label{thm:main}
Keep the notation above. 
The VOA $\Vam$ is simple if and only if 
$r \in \BC$ is not an integer, that is, 
$r \in \BC \stm \BZ$.
\end{thm}

%
We will prove Theorem~\ref{thm:main} as follows. 
First, in \S\ref{sec:step1}, 
we will show that $\Vam \subset \Mam$ is, in fact, 
identical to the whole of $\Mam$ (Proposition~\ref{V=M}),
and then prove that the VOA $\Vam$ ($=\Mam$) is simple 
if and only if $\Mam$ is irreducible as an $\Lr$-module 
(Proposition~\ref{prop:step1}). 
Let $\Lon$ be a Lie subalgebra of $\Lr$ 
generated by $\bigl\{v^{11}(m,\,n) \mid 
m,\,n \in \BZ \text{ with } m \le n\bigr\} \subset \CB$, and 
set $\Mon=U(\Lon)\vm \subset \Mam$. 
In \S\ref{sec:step2}, it will be shown that 
the $\Lr$-module $\Mam$ is irreducible
if and only if $\Mon$ is 
irreducible as an $\Lon$-module (Proposition~\ref{prop:step2}). 
In \S\ref{sec:step3}, we will prove that 
if $r \in \BC \stm \BZ$, then 
$\Mon$ is an irreducible $\Lon$-module, and 
hence $\Vam$ is a simple VOA (Proposition~\ref{prop:step3}). 
Finally, in \S\ref{sec:step4}, we will give 
some singular vectors of the $\Lon$-module $\Mon$ 
explicitly in the case that $r \in \BZ$ (Proposition~\ref{prop:step4}), 
which implies that $\Mon$ is reducible, and hence $\Vam$ is not simple. 
%
%
\section{Simplicity of $\Vam$ and irreducibility of $\Mam$.}
\label{sec:step1}
%
%
\subsection{Relation between $\Vam$ and $\Mam$.}
\label{subsec:V=M}
As in the previous section, 
we fix $d \in \BZ_{\ge 2}$ and $r \in \BC$. 
This subsection is devoted to proving 
the following proposition. 
%
%
\begin{prn} \label{V=M}
The subspace $\Vam \subset \Mam$ is identical to 
the whole of $\Mam$, that is, $\Vam=\Mam$ holds. 
\end{prn}


In order to prove Proposition~\ref{V=M},
we need some technical lemmas.
%
%
\begin{lem}\label{lem3:1}
{\rm (1)} 
For each $1 \le i,\,j \le d$, we have 
%
%
\begin{equation} \label{vij1}
v^{ij}(-1,\,-1)\vm = 2\lij(-2)\vm.
\end{equation}

\noindent
{\rm (2)}
Let $1 \le i,\,j \le d$ with $i \ne j$, and 
$m,\,n \in \BZ_{< 0}$. Then, 
%
%
\begin{equation} \label{vij}
v^{ij}(m-1,\,n) \vm =-\frac{1}{m}\li(-1)v^{ij}(m,\,n)\vm.
\end{equation}

\noindent
{\rm (3)} 
Let $1 \le i,\,j \le d$ with $i \ne j$, and 
$m,\,n \in \BZ_{< 0}$. Then, 
%
%
\begin{equation} \label{vii}
v^{ii}(m-1,\,n) \vm
=\frac{2}{m(m-n+1)}\li(0)\lij(-1)v^{ij}(n,\,m)\vm.
\end{equation}
\end{lem}

\begin{proof}
(1) By the definition \eqref{dLij} of $\lij(m)$,
\begin{equation*}
2\lij(-2)\vm=\sum_{h \in \BZ}v^{ij}(-2-h,\,h) \vm.
\end{equation*}
Note that $v^{ij}(-2-h,\,h)=v^{ji}(h,\,-2-h)$ for all $h \in \BZ$ 
(see Remark~\ref{rem:CB}). 
Since $v^{ij}(m,\,n)\vm=0$ if $v^{ij}(m,\,n) \in \CB_{+}$, 
it follows that $v^{ij}(-2-h,\,h)\vm=0$ unless $h=-1$.
Thus we obtain $2\lij(-2)\vm=v^{ij}(-1,\,-1)\vm$, 
and hence \eqref{vij1}. 

\noindent
(2) As in the proof of part (1), we can easily show that
%
%
\begin{equation} \label{eq:lij-1}
\lij(-1)\vm
=\frac{1}{2}\sum_{h \in \BZ}v^{ij}(-1-h,\,h)\vm=0
\end{equation}
for all $1 \le i,\,j \le d$. 
Hence, 
\begin{align*}
\li(-1)v^{ij}(m,\,n)\vm 
& = [\li(-1),\, v^{ij}(m,\,n)] \vm
  \qquad 
  \text{since $\li(-1)\vm=0$ by \eqref{eq:lij-1}} \\[1.5mm]
& =\frac{1}{2}\sum_{h \in \BZ}
   [v^{ii}(-1-h,\,h),\,v^{ij}(m,\,n)] \vm.
\end{align*}
By direct computation (as in Example~\ref{eg:bracket}), 
we see that
\begin{align*}
& [v^{ii}(-1-h,\,h),\,v^{ij}(m,\,n)]_{r} = \\
& \hspace*{20mm}
\delta_{-1-h+m,0}\,(-1-h) v^{ij}(h,\,n) + 
\delta_{h+m,0}\,h v^{ij}(-1-h,\,n). 
\end{align*}
Therefore, 
\begin{align*}
& \frac{1}{2}\sum_{h \in \BZ}
   [v^{ii}(-1-h,\,h),\,v^{ij}(m,\,n)] \vm \\[1.5mm]
& \qquad =
  \frac{1}{2}\sum_{h \in \BZ}
  \bigl\{
    \delta_{-1-h+m,0}\,(-1-h) v^{ij}(h,\,n) + 
    \delta_{h+m,0}\,h v^{ij}(-1-h,\,n)
  \bigr\}\vm \\[1.5mm]
& \qquad = 
  \frac{1}{2}
  \bigl\{(-m) v^{ij}(m-1,\,n) + (-m) v^{ij}(m-1,\,n) \bigr\}\vm \\[1.5mm]
& \qquad = 
  (-m) v^{ij}(m-1,\,n)\vm.
\end{align*}
Thus we obtain $\li(-1)v^{ij}(m,\,n)\vm = (-m) v^{ij}(m-1,\,n)$, 
and hence \eqref{vij}. 

\noindent
(3) We have
\begin{align*}
\lij(-1)v^{ij}(n,\,m)\vm & =
  [\lij(-1),\,v^{ij}(n,\,m)]\vm
  \qquad \text{since $\lij(-1)\vm=0$ by \eqref{eq:lij-1}} \\[1.5mm]
& =\frac{1}{2} \sum_{h \in \BZ}
  [v^{ij}(-1-h,\,h),\,v^{ij}(n,\,m)] \vm.
\end{align*}
As in Example~\ref{eg:bracket}, 
we see that
\begin{align*}
& [v^{ij}(-1-h,\,h),\,v^{ij}(n,\,m)]_{r} = \\
& \hspace*{20mm}
\delta_{h+m,0}\,h v^{ii}(n,\,-1-h) + 
\delta_{-1-h+n,0}\,(-1-h) v^{jj}(h,\,m). 
\end{align*}
Therefore we get
\begin{align*}
& \frac{1}{2} \sum_{h \in \BZ}
  [v^{ij}(-1-h,\,h),\,v^{ij}(n,\,m)] \vm \\[1.5mm]
& \qquad = 
  \frac{1}{2}\sum_{h \in \BZ}
  \bigl\{
    \delta_{h+m,0}\,h v^{ii}(n,\,-1-h) + 
    \delta_{-1-h+n,0}\,(-1-h) v^{jj}(h,\,m)
  \bigr\}\vm \\[1.5mm]
& \qquad =
  \frac{1}{2} 
  \bigl\{
    (-m) v^{ii}(n,\,m-1)\vm + 
    (-n) v^{jj}(n-1,\,m)\vm
  \bigr\}, 
\end{align*}
and hence 
\begin{equation*}
\lij(-1)v^{ij}(n,\,m)\vm = 
  -\frac{m}{2} v^{ii}(n,\,m-1)\vm 
  -\frac{n}{2} v^{jj}(n-1,\,m)\vm.
\end{equation*}
Now, since $\li(0)\vm=0$, it follows that 
\begin{align*}
& \li(0)\lij(-1)v^{ij}(n,\,m)\vm \\[2mm]
& \qquad =
  -\frac{m}{2} \li(0) v^{ii}(n,\,m-1)\vm 
  -\frac{n}{2} \li(0) v^{jj}(n-1,\,m)\vm \\[2mm]
& \qquad = 
  -\frac{m}{2} [\li(0),\,v^{ii}(n,\,m-1)]\vm 
  -\frac{n}{2} [\li(0),\,v^{jj}(n-1,\,m)]\vm.
\end{align*} 
It can be easily checked that 
$[\li(0),\,v^{jj}(n-1,\,m)]=0$ since $i \ne j$. 
Thus, 
\begin{align*}
& -\frac{m}{2} [\li(0),\,v^{ii}(n,\,m-1)]\vm 
  -\frac{n}{2} [\li(0),\,v^{jj}(n-1,\,m)]\vm \\[2mm]
& \hspace*{5mm} 
  = -\frac{m}{2} [\li(0),\,v^{ii}(n,\,m-1)]\vm \\[2mm]
& \hspace*{5mm}
  = \biggl\{
      -\frac{m}{4} \underbrace{[v^{ii}(0,\,0),\,v^{ii}(n,\,m-1)]}_{=0}
      -\frac{m}{2} \sum_{h \in \BZ_{> 0}}
      [v^{ii}(-h,\,h),\,v^{ii}(n,\,m-1)] \biggr\}\vm. 
\end{align*}
Since $n,\,m-1 \in \BZ_{< 0}$, it follows that 
\begin{equation*}
[v^{ii}(-h,\,h),\,v^{ii}(n,\,m-1)]_{r} = 
\delta_{h+n,0}\,h v^{ii}(-h,\,m-1) + 
\delta_{h+m-1,0}\,h v^{ii}(n,\,-h)
\end{equation*}
for $h \in \BZ_{\ge 1}$. Therefore, 
\begin{align*}
& -\frac{m}{2} 
  \sum_{h \in \BZ_{> 0}}
  [v^{ii}(-h,\,h),\,v^{ii}(n,\,m-1)] \vm \\[2mm]
& \qquad = 
  -\frac{m}{2} 
  \sum_{h \in \BZ_{> 0}}
  \bigl\{
  \delta_{h+n,0}\,h v^{ii}(-h,\,m-1) + 
  \delta_{h+m-1,0}\,h v^{ii}(n,\,-h)
  \bigr\}\vm \\[2mm]
& \qquad = 
  -\frac{m}{2} 
  \bigl\{
   (-n) v^{ii}(n,\,m-1) + 
   (-m+1) v^{ii}(n,\,m-1)
  \bigr\}\vm \\[2mm]
& \qquad = 
  \frac{m(m+n-1)}{2} v^{ii}(m-1,\,n) \vm.
\end{align*}
Thus we obtain 
\begin{equation*}
\li(0)\lij(-1)v^{ij}(n,\,m)\vm = 
\frac{m(m+n-1)}{2} v^{ii}(m-1,\,n) \vm, 
\end{equation*}
and hence \eqref{vii}. 
This completes the proof of the lemma.
\end{proof}
%
%
\begin{lem} \label{lem3:4}
{\rm (1)} 
Let $1 \le i,\,j \le d$ with $i \ne j$, and 
$m,\,n \in \BZ_{< 0}$. Then, 
%
%
\begin{equation} \label{vjl}
v^{ij}(m,\,n)\vm= \alpha 
\li(-1)^{-m-1} \lj(-1)^{-n-1} \lij(-2)\vm
\end{equation}
for some $\alpha \in \BC \stm \{0\}$.

\noindent
{\rm (2)} 
Let $1 \le i \le d$, and $m,\,n \in \BZ_{< 0}$ with $m \le -2$. 
Take $1 \le j \le d$ with $j \ne i$ arbitrarily. 
Then, 
%
%
\begin{equation} \label{vil}
v^{ii}(m,\,n)\vm = \beta 
\li(0)\lij(-1) \li(-1)^{-n-1} \lj(-1)^{-m-2}\lij(-2)\vm
\end{equation}
for some $\beta \in \BC \stm \{0\}$. 
\end{lem}

\begin{proof}
(1) We prove \eqref{vjl} 
by induction on $-m-n$ (note that $-m-n \geq 2$).
If $-m-n=2$, that is, $m=n=-1$,
then \eqref{vjl} follows immediately from \eqref{vij1}.
Assume that $-m-n>2$. By Remark~\ref{rem:CB}, 
we may assume that $m < -1$. 
Then, by \eqref{vij}, we have
\begin{equation*}
v^{ij}(m,\,n)\vm = 
 - \frac{1}{m+1}\li(-1)v^{ij}(m+1,\,n)\vm.
\end{equation*}
Applying the inductive assumption to 
the right-hand side of the equation above,
we obtain
\begin{align*}
v^{ij}(m,\,n)\vm
 & = -\frac{1}{m+1}\li(-1)v^{ij}(m+1,\,n)\vm \\[2mm]
 & = -\frac{1}{m+1}\li(-1)
     \left\{ 
       \alpha \li(-1)^{-m-2} \lj(-1)^{-n-1} \lij(-2)\vm
     \right\} \\[2mm]
 & = - \frac{\alpha}{m+1} \li(-1)^{-m-1} \lj(-1)^{-n-1}\lij(-2)\vm,
\end{align*}
where $\alpha \in \BC \stm \{0\}$. 
Thus we have proved part (1). 

\noindent 
(2) Using \eqref{vii} and \eqref{vjl}, we have
\begin{align*}
& v^{ii}(m,\,n)\vm = 
  \frac{2}{(m+1)(m+n)}\li(0)\lij(-1)v^{ij}(n,\,m+1)\vm 
  \qquad \text{by \eqref{vii}} \\[2mm]
& = \frac{2\beta}{(m+1)(m+n)}
  \li(0)\lij(-1)\li(-1)^{-n-1}\lj(-1)^{-m-2}\lij(-2)\vm
  \qquad  \text{by \eqref{vil}},
\end{align*}
where $\beta \in \BC \stm \{0\}$. 
Thus we have proved part (2), thereby completing 
the proof of Lemma~\ref{lem3:4}.
\end{proof}

\begin{proof}[Proof of Proposition~\ref{V=M}]
We set
\begin{equation*}
U:=\Span_{\BC} 
  \left\{
    \lab{i_1j_1}{m_1}\lab{i_2j_2}{m_2} \cdots \lab{i_pj_p}{m_p}\vm 
    \ \Biggm| 
    \begin{array}{l}
    p \geq 0, \text{ and } 1 \leq i_q,\,j_q \leq d, \\[2mm]
    m_q \in \{-2,-1,0\} \text{ for } 1 \leq q \leq p
    \end{array}
  \right\}.
\end{equation*}
In order to prove Proposition~\ref{V=M},
it suffices to show that $U=\Mam$.
Indeed, it is obvious from the definition of $\Vam$ that 
$U \subset \Vam$. Hence, if $U=\Mam$ holds, then 
$\Mam=U \subset \Vam \subset \Mam$,
which implies that $\Vam=\Mam$.

\begin{claim*}
We have $x U \subset U$ 
for every $x \in \CB$.
\end{claim*}

\noindent 
{\it Proof of Claim.} 
Fix $x \in \CB$. It suffices to prove that
\begin{equation*}
x \lab{i_1j_1}{m_1}\lab{i_2j_2}{m_2}\cdots
\lab{i_pj_p}{m_p}\vm \in U
\end{equation*}
for every $p \geq 0$ and every $1 \leq i_q,\,j_q \leq d$,
$m_q \in \{-2, -1, 0\}$ for $1 \leq q \leq p$. 
We prove this by induction on $p$. 
It is obvious that $x\vm=0 \in U$ if $x \in \CB_{+}$. 
Also, we see from equations \eqref{vjl} and \eqref{vil}
that $x\vm$ is contained in $U$ if $x \in \CB_{-}$. 
Thus the assertion holds when $p=0$.
Assume that $p > 0$. Then we have
%
%
\begin{align}
& x \lab{i_1j_1}{m_1}\lab{i_2j_2}{m_2} \cdots 
  \lab{i_pj_p}{m_p}\vm = \notag \\
& [x,\,\lab{i_1j_1}{m_1}]
    \lab{i_2j_2}{m_2}\cdots \lab{i_pj_p}{m_p}\vm 
  +\lab{i_1j_1}{m_1} 
   \bigl\{ x \lab{i_2j_2}{m_2}\cdots \lab{i_pj_p}{m_p}\vm \bigr\}.
  \label{indf}
\end{align}
Since $x \lab{i_2j_2}{m_2}\cdots \lab{i_pj_p}{m_p}\vm$ is 
contained in $U$ by the inductive assumption, and 
since $m_1 \in \{-2,-1,0\}$, 
it follows that the second term of 
the right-hand side of \eqref{indf} is contained in $U$.
Now we deduce from the definition of 
the Lie bracket on $\CL_{r}$ and the 
definition \eqref{dLij} of $L^{ij}_{r}(m)$ that 
$[x,\,\lab{i_1j_1}{m_{1}}]$ 
can be written in the form: 
\begin{equation*}
[x,\,\lab{i_1j_1}{m_1}]
= \alpha_{1} y_{1} + \alpha_{2}y_{2} + \cdots + \alpha_{s}y_{s} + \beta
\end{equation*}
for some $\alpha_{1},\,\alpha_{2},\,\dots,\,\alpha_{s},\,\beta \in \BC$ 
and $y_{1},\,y_{2},\,\dots,\,y_{s} \in \CB$. 
By substituting this into the first term of 
the right-hand side of \eqref{indf}, we see that
\begin{align*}
& [x,\,\lab{i_1j_1}{m_1}]
  \lab{i_2j_2}{m_2} \cdots \lab{i_pj_p}{m_p}\vm \\[2mm]
& = 
  \sum_{1 \le s' \le s}\alpha_{s'} 
  \underbrace{
    y_{s'}\lab{i_2j_2}{m_2}\cdots \lab{i_pj_p}{m_p}\vm
    }_{\in U \ \text{by the inductive assumption}}+
  \beta
  \underbrace{
    \lab{i_2j_2}{m_2}\cdots \lab{i_pj_p}{m_p}\vm
    }_{\in U}, 
\end{align*}
and hence that the first term also is contained in $U$.
Therefore we conclude that the left-hand side of \eqref{indf} 
is contained in $U$, thereby completing the proof of Claim.

\vspace{3mm}

The claim above implies that 
$U$ is an $\CL_{r}$-submodule of $\Mam$ which contains $\vm$.
Hence we conclude that $U=\Mam$, thereby completing 
the proof of Proposition~\ref{V=M}.
\end{proof}
%
%
\subsection{Relation between the simplicity of $\Vam$
and the irreducibility of $\Mam$.}
\label{subsec:Vs-Mir}
In this subsection, we prove the following proposition. 
%
%
\begin{prn} \label{prop:step1}
The VOA $\Vam \ (= \Mam)$ is simple if and only if 
$\Mam$ is irreducible as an $\Lr$-module.
\end{prn}

First let us show the following lemma,
needed in the proof of Proposition~\ref{prop:step1}. 
%
%
\begin{lem}\label{lem3:5}
Let $1 \le i,\,j \le d$ with $i \ne j$, and 
$m,\,n \in \BZ_{< 0}$. The vertex operator 
$Y(v^{ij}(m,\,n)\vm, z)$ of 
$v^{ij}(m,\,n)\vm \in \Mam=\Vam$ is given by\,{\rm:}
%
%
\begin{align}
& Y(v^{ij}(m,\,n)\vm, z) = 
   (-1)^{-m-n} \times \notag \\[2mm]
& \sum_{l \in \BZ}
  \left\{
   \sum_{k \in \BZ}
      \binom{l+n-k}{-m-1}
      \binom{k-n-1}{-n-1}
   v^{ij}(l+m+n+1-k,\,k)
  \right\} z^{-l-1}. \label{voij}
\end{align}
\end{lem}

\begin{proof}
We prove the lemma by induction on $-m-n$
(note that $-m-n \geq 2$).
If $-m-n=2$, that is, $m=n=-1$,
then it follows that
\begin{align*}
& Y(v^{ij}(-1,\,-1)\vm, z) = 
  2 Y(L^{ij}_{r}(-2)\vm, z) \quad \text{by \eqref{vij1}} \\[2mm]
& \hspace*{15mm}
  = 2 Y(\omega^{ij}_{r}, z) 
  = 2 \sum_{l \in \BZ}\lij(l)z^{-l-2} 
   \qquad \text{by \eqref{eq:y0} and Theorem~\ref{thm:am}} \\[2mm]
& \hspace*{15mm}
   = \sum_{l \in \BZ}
     \left\{
       \sum_{k \in \BZ}v^{ij}(l-k,\,k)
     \right\}z^{-l-2}
   = \sum_{l \in \BZ}
     \left\{
       \sum_{k \in \BZ}v^{ij}(l-1-k,\,k)
     \right\}z^{-l-1}.
\end{align*}
Therefore the equation \eqref{voij} holds if $-m-n=2$.
Assume that $-m-n>2$. By Remark~\ref{rem:CB}, 
we may assume that $m < -1$.
It follows from \eqref{vij} that 
\begin{equation*}
Y(v^{ij}(m,\,n)\vm, z) = 
  \frac{1}{-m-1}Y(\li(-1)v^{ij}(m+1,\,n)\vm, z).
\end{equation*}
By using the commutator formula 
(see \cite[p.\,54]{ll}), we deduce that 
\begin{align*}
(\li(-1)v^{ij}(m+1,\,n)\vm)_{l} 
 & = ((\oi)_{0}v^{ij}(m+1,\,n)\vm)_l 
   = [(\oi)_{0},\,(v^{ij}(m+1,\,n)\vm)_l] \\
 & = [\li(-1),\,(v^{ij}(m+1,\,n)\vm)_l].
\end{align*}
Also, it follows from the inductive assumption that 
\begin{align*}
& (v^{ij}(m+1,\,n)\vm)_l 
  = (-1)^{-m-1-n} \times \notag \\[2mm]
& \hspace*{20mm}
  \left\{
   \sum_{k \in \BZ}
      \binom{l+n-k}{-m-2}
      \binom{k-n-1}{-n-1}
   v^{ij}(l+m+n+2-k,\,k)
  \right\}.
\end{align*}
Combining these equations, we get 
\begin{align*}
& Y(v^{ij}(m,\,n)\vm, z) = \frac{(-1)^{-m-1-n}}{-m-1}
  \sum_{l \in \BZ} \, 
  \Biggl\{ \sum_{k \in \BZ}
      \binom{l+n-k}{-m-2}
      \binom{k-n-1}{-n-1} \times \\[2mm]
& \hspace*{40mm}
  \bigl[\li(-1),\,v^{ij}(l+m+n+2-k,\,k)\bigr]\Biggr\} z^{-l-1}.
\end{align*}
It can be easily checked by direct computation that 
\begin{align*}
& \bigl[\li(-1),\,v^{ij}(l+m+n+2-k,\,k)\bigr] \\[2mm]
& \hspace*{10mm} =
  \frac{1}{2} 
  \sum_{h \in \BZ} 
  \bigl[v^{ii}(-1-h,\,h),\,v^{ij}(l+m+n+2-k,\,k)\bigr] \\[2mm]
& \hspace*{10mm} = 
  -(l+m+n+2-k)\,v^{ij}(l+m+n+1-k,\,k).
\end{align*}
Thus we obtain 
\begin{align*}
& Y(v^{ij}(m,\,n)\vm, z) \\[2mm]
& = (-1)^{-m-n}
  \sum_{l \in \BZ}\, 
  \Biggl\{ \sum_{k \in \BZ}
      \frac{1}{-m-1}\binom{l+n-k}{-m-2}(l+m+n+2-k) \times \\[2mm]
& \hspace*{40mm}
      \binom{k-n-1}{-n-1}
      v^{ij}(l+m+n+1-k,\,k)\Biggr\} z^{-l-1} \\[2mm]
& = (-1)^{-m-n}
  \sum_{l \in \BZ}\, 
  \Biggl\{ \sum_{k \in \BZ}
      \binom{l+n-k}{-m-1}
      \binom{k-n-1}{-n-1}
      v^{ij}(l+m+n+1-k,\,k)
  \Biggr\} z^{-l-1}. 
\end{align*}
Thus we have proved Lemma~\ref{lem3:5}. 
\end{proof}

\begin{proof}[Proof of Proposition~\ref{prop:step1}]
First, we show the ``if'' part.
Assume that $\Vam=\Mam$ is not simple,
and let $W \subset \Vam=\Mam$ be a proper ideal of the VOA $\Vam$.
Let us show that $W$ is a (proper) 
$\Lr$-submodule of $\Mam$.
%
%
\begin{claim} \label{c:s1-1}
Let $1 \le i,\,j \le d$ with $i \ne j$. 
Then, $v^{ij}(m,\,n)W \subset W$ for all $m,\,n \in \BZ$.
\end{claim}

\noindent
{\it Proof of Claim~\ref{c:s1-1}.} Let $u \in W$. 
We prove that $v^{ij}(m,\,n)u \in W$ for all $m, n \in \BZ$, 
which is equivalent to showing that 
$v^{ij}(s-t,\,t)u \in W$ for all $s, t \in \BZ$.
Fix $s \in \BZ$. If $v^{ij}(s-t,\,t)u=0$ for all $t \in \BZ$, 
then the assertion is obvious. So, let us assume that 
$v^{ij}(s-t,\,t)u \neq 0$ for some $t \in \BZ$. 
Take $t_1,\,t_2 \in \BZ$ with $t_1 \leq t_2$ in such a way
that if $v^{ij}(s-t,\,t)u \neq 0$, then $t_1 \leq t \leq t_2$.
By Lemma \ref{lem3:5}, we see that 
\begin{align*}
\bigl(v^{ij}(-p,\,-1)\vm\bigr)_{s+p}u
 & = (-1)^{p+1}\sum_{t \in \BZ}
     \binom{s+p-1-t}{p-1}
     v^{ij}(s-t,\,t)u \\[2mm]
 & = \sum_{t_{1} \le t \le t_{2}} (-1)^{p+1}
     \binom{s+p-1-t}{p-1}
     v^{ij}(s-t,\,t)u. 
\end{align*}
for $p>0$. Because $W$ is an ideal of the VOA $\Vam=\Mam$, 
it follows that 
\begin{equation*}
\bigl(v^{ij}(-p,\,-1)\vm\bigr)_{s+p}u \in W
\end{equation*}
for all $p>0$. Hence, in order to prove that 
$v^{ij}(s-t,\,t)u \in W$ for all $t_1 \leq t \leq t_2$, 
it suffices to show that the matrix
%
%
\begin{equation}\label{mat}
\left((-1)^{p+1}
\binom{s+p-1-t}{p-1}\right)_{
  \begin{subarray}{c}
  1 \leq p \leq t_2-t_1+1 \\[1mm]
  t_1 \leq t \leq t_2
  \end{subarray}}
\end{equation}
is invertible. So we show the following: 
Let $L \in \BZ$, and $M \in \BZ_{> 0}$. Then, 
%
%
\begin{equation} \label{eq:det}
\det(a_{p,N})_{1 \leq p,N \leq M} \neq 0, \quad
\text{where } 
a_{p,N}:=\binom{L+p-N}{p-1}.
\end{equation}
We prove \eqref{eq:det} by induction on the size $M$ of 
the matrix. The claim is obvious when $M=1$.
Assume that $M>1$. Using the formula 
\begin{equation*}
a_{p,N}-a_{p,N-1}=-
\binom{L+p-N}{p-2}
\quad \text{for $2 \leq N \leq M$},
\end{equation*}
we deduce that
\begin{align*}
\det(a_{p,N})_{1\leq p,\,N \leq M}
& = \det \left(-\binom{L+p-N}{p-2}\right)_{2 \leq p,N \leq M} \\[2mm]
& = (-1)^{M-1} \det \left(\binom{L+p-N}{p-2}\right)_{2 \leq p,N \leq M} \\[2mm]
& = (-1)^{M-1} 
    \underbrace{\det(a_{p',N'})_{1 \leq p',N' \leq M-1}}_{
    \neq 0 \, \text{by the inductive assumption}},
\end{align*}
where $p':=p-1$ and $N':=N-1$. 
Thus we obtain \eqref{eq:det}.
We see from \eqref{eq:det} with 
$L=s-t_{1}$, $M=t_2-t_1+1$, and $N=t-t_1+1$ that 
the determinant of the matrix \eqref{mat}
is not equal to $0$,
and hence the matrix (\ref{mat}) 
is invertible. Thus we have proved Claim~\ref{c:s1-1}. 

\vspace{3mm}

%
\begin{claim} \label{c:s1-2}
Let $1 \le i \le d$. Then, 
$v^{ii}(m,\,n)W \subset W$ for all $m,n \in \BZ$.
\end{claim}

\noindent
{\it Proof of Claim~\ref{c:s1-2}.} Let $u \in W$. 
Take $1 \le j \le d$ with $j \ne i$ arbitrarily, and take 
$N \in \BZ_{<0}$ in such a way that $v^{jj}(N,\,-N)u=0$. 
By direct computation, we see that 
\begin{align*}
[v^{ij}(m,\,-N),\,v^{ij}(n,\,N)] u
 & = -N v^{ii}(m,\,n)u+\delta_{m+n,0}mv^{jj}(N,\,-N)u+\alpha u \\
 & = -N v^{ii}(m,\,n)u + \alpha u \quad 
   \text{since $v^{jj}(N,\,-N)u=0$}
\end{align*}
for some $\alpha \in \BC$.
Since $[v^{ij}(m,\,-N),\,v^{ij}(n,\,N)]u$ is contained in $W$ 
by Claim~\ref{c:s1-1}, we conclude that 
$v^{ii}(m,\,n)u \in W$, thereby 
completing the proof of Claim~\ref{c:s1-2}.

\vspace{3mm}

It follows from Claims~\ref{c:s1-1} and \ref{c:s1-2} that 
$v^{ij}(m,\,n)W \subset W$ for all $1 \le i,\,j \le d$ and 
$m,n \in \BZ$, which implies that $W$ is an $\CL_{r}$-submodule 
of $M_{r}$. Thus we have proved the ``if'' part of 
Proposition~\ref{prop:step1}. 

Next, we show the ``only if'' part.
Assume that $\Mam=\Vam$ is a reducible $\Lr$-module,
and let $W \subset \Mam=\Vam$ be a proper $\Lr$-submodule.
Let $v \in \Vam$, and $l \in \BZ$. 
By the definition of the vertex operator of $\Vam$,
we deduce that $v_l \in \End(\Vam)$ can be 
written as an (infinite) linear combination of 
products of $v^{ij}(m,\,n)$, $1 \le i,\,j \le d$, $m,\,n \in \BZ$.
Since $v^{ij}(m,\,n)W \subset W$ by assumption, 
we have $v_l W \subset W$. Therefore, 
$W$ is a proper ideal of the VOA $\Vam$,
and hence $\Vam$ is not simple. 
This completes the proof of Proposition~\ref{prop:step1}. 
\end{proof}

The next corollary follows immediately from 
the proof of Proposition~\ref{prop:step1}.
%
%
\begin{cor} \label{cor:max1}
Assume that the VOA $\Vam$ is not simple, or equivalently, 
the $\CL_{r}$-module $\Mam$ is reducible. 
Then, $W \subset \Vam \ (=M_{r})$ is 
the maximal proper ideal of the VOA $\Vam$ if and only if 
$W$ is the maximal proper $\CL_{r}$-submodule of $\Mam$. 
\end{cor}

%
\section{Irreducibility of $\Mam$ and $\Mon$.}
\label{sec:step2}

Let $\Lon$ be the Lie subalgebra of $\Lr$ generated by 
$\CB^{(1)}:=
 \bigl\{v^{11}(m,n) \mid 
 m,\,n \in \BZ \text{ with } m \le n \bigr\} \subset \CB$. 
Set $\Mon=U(\Lon)\vm \subset \Mam$. In this section,
we prove the following proposition.
%
%
\begin{prn}\label{prop:step2}
The $\Lr$-module $\Mam$ is irreducible
if and only if $\Mon$ is irreducible as 
an $\Lon$-module.
\end{prn}
%
%
\subsection{Proof of the ``if'' part of 
Proposition~\ref{prop:step2}.}
\label{subsec:s2-if}

In this subsection,
we show that if $\Mon$ is an irreducible $\Lon$-module,
then $\Mam$ is an irreducible $\Lr$-module.
This assertion follows immediately from the next lemma.
%
%
\begin{lem}\label{lem4:2}
Let $W$ be a nonzero $\Lr$-submodule of $\Mam$.
Then, $W \cap \Mon$ is a nonzero $\Lon$-submodule of $\Mon$.
\end{lem}

Indeed, if $W \subset \Mam$ is a nonzero $\Lr$-submodule of $\Mam$, 
then it follows from Lemma~\ref{lem4:2} that 
$W \cap \Mon$ is a nonzero $\Lon$-submodule of $\Mon$. 
Since $\Mon$ is assumed to be an irreducible $\Lon$-module, 
we have $W \cap \Mon=\Mon$. In particular, 
$W$ contains $\vm \in \Mam$, which implies that $W=\Mam$.


In order to prove Lemma~\ref{lem4:2}, we introduce 
a weight space decomposition of $\Mam$. Define 
\begin{equation*}
\CH:=
 \bigoplus_{1 \le k \le d,\, l \in \BZ_{< 0}} 
 \BC v^{kk}(l,\,-l) \subset \Lr.
\end{equation*}
Then it can be easily seen 
that $\CH$ is an abelian Lie subalgebra of $\Lr$. 
Set $h_{k,\,l}:=-(1/l)\,v^{kk}(l,\,-l) \in \CH$ 
for $1 \le k \le d$ and $l \in \BZ_{< 0}$. 
By simple computation, we see that 
for $1 \le k \le d$ and $l \in \BZ_{< 0}$, 
$1 \le i,\,j \le d$ and $m,\,n \in \BZ$, 
%
%
\begin{equation} \label{eq:hkl}
[h_{k,\,l},\,v^{ij}(m,\,n)]_{r}= 
 \bigl(
   \delta_{(k,\,l),\,(i,m)}+\delta_{(k,\,l),\,(j,n)}-
   \delta_{(k,-l),\,(i,m)}-\delta_{(k,-l),\,(j,n)}
 \bigr) v^{ij}(m,\,n).
\end{equation}
Let $\Lambda_{k,\,l} \in \CH^{\ast}:=\Hom_{\BC}(\CH,\,\BC)$ 
be the dual basis of $h_{k,\,l} \in \CH$ 
for $1 \le k \le d$ and $l \in \BZ_{< 0}$, and set
\begin{equation*}
Q_{+}:=
  \sum_{1 \le k \le d,\, l \in \BZ_{< 0}}
  \BZ_{\ge 0} \Lambda_{k,\,l} \subset \CH^{\ast}. 
\end{equation*}
We see from \eqref{eq:hkl} that 
for each $\bx \in \CS$, the basis element $w(\bx)$ of $\Mam$ 
is contained in the ``weight space'' 
$(M_{r})^{\lambda}:=\bigl\{u \in M_{r} \mid 
 hu=\lambda(h)u \text{ for all } h \in \CH \bigr\}$ 
of weight $\lambda$ for some $\lambda \in Q_{+}$ 
(see also Remark~\ref{rem:wt2} below). 
Thus the $\Lr$-module $M_{r}$ admits 
the weight space decomposition with respect to 
the abelian Lie subalgebra $\CH \subset \Lr$ as:
%
%
\begin{equation} \label{eq:wsd}
M_{r}=\bigoplus_{\lambda \in Q_{+}} (M_{r})^{\lambda}.
\end{equation}
%
%
\begin{rmk} \label{rem:wt}
Let $1 \le i,\,j \le d$, and 
let $m,\,n \in \BZ$. 
It follows from \eqref{eq:hkl} that 
$v^{ij}(m,n) \in \End_{\BC}(M_{r})$ is a 
homogeneous operator of weight 
$-\sn(m) (\Lambda_{i,-m}+\Lambda_{i,m}) 
 -\sn(n) (\Lambda_{j,-n}+\Lambda_{j,n})$, 
where for $N \in \BZ$, 
\begin{equation*}
\sn(N):=
 \begin{cases}
 1 & \text{if $N > 0$}, \\
 0 & \text{if $N=0$}, \\
 -1 & \text{if $N < 0$}, 
 \end{cases}
\end{equation*}
and for convenience, 
$\Lambda_{k,l}:=0$ for all $1 \le k \le d$ 
and $l \in \BZ_{\ge 0}$. 
\end{rmk}
%
%
\begin{rmk} \label{rem:wt2}
Let $\bx=(x_{p} \succeq x_{p-1} \succeq \cdots \succeq x_1) \in \CS$ 
with $x_q=v^{i_{q}j_{q}}(m_{q},\,n_{q}) \in \CB_{-}$ for $1 \le q \le p$. 
For $1 \le k \le d$ and $l \in \BZ_{< 0}$, we define 
\begin{align*} 
& \nu(\bx,\,(k,l)):= \\
& \quad 
\# \bigl\{1 \le q \le p \mid (i_q,\,m_q)=(k,\,l)\bigr\}
 +\# \bigl\{1 \le q \le p \mid (j_q,\,n_q)=(k,\,l)\bigr\}.
\end{align*}
Namely, $\nu(\bx,\,(k,l))$ denotes ``the number of 
$v^{k}(l)$ appearing in $\bx$'' (see \eqref{eq:dvij}). 
Then we deduce from Remark~\ref{rem:wt} that 
the weight of $w(\bx)$ is equal to 
\begin{equation*}
  \sum_{1 \le k \le d,\, l \in \BZ_{< 0}}
  \nu(\bx,\,(k,l)) \Lambda_{k,\,l} \in Q_{+}. 
\end{equation*}
\end{rmk}

\begin{proof}[Proof of Lemma~\ref{lem4:2}] 
Remark that the submodule $W \subset \Mam$ also admits 
the weight space decomposition 
$W=\bigoplus_{\lambda \in Q_{+}} W^{\lambda}$, 
where $W^{\lambda}:=W \cap (\Mam)^{\lambda}$. 
Let $u$ be a nonzero homogeneous element of $W$, 
that is, $u \in W^{\lambda} \stm \{0\}$ 
for some $\lambda \in Q_{+}$.
It suffices to show that there exists $x \in U(\CL_{r})$ 
such that $xu \in \Mon \stm \{0\}$. 
For each $\xi \in Q_{+}$, we set 
%
%
\begin{equation} \label{eq:theta}
\theta(\xi):=\sum_{2 \le k \le d,\, l \in \BZ_{< 0}} 
  \xi(h_{k,\,l}).
\end{equation}
We show the claim above by induction on $\theta(\lambda)$. 
If $\theta(\lambda) = 0$, then the claim is 
obvious since $u \in M_{r}^{(1)}$ (see Remark~\ref{rem:wt2}). 
Assume that $\theta(\lambda) > 0$. 
Let $2 \le i \le d$ and $m \in \BZ_{< 0}$ 
be such that $\lambda(h_{i,m}) \ge 1$, and let 
$N \in \BZ_{< 0}$ be such that $\lambda(h_{1,N})=0$. 
Then we have $v^{1i}(-N,\,m)v^{1i}(N,\,-m)u \ne 0$. 
Indeed, since 
\begin{equation*}
\bigl[v^{1i}(-N,\,m),\,v^{1i}(N,\,-m)\bigr]_{r} =
  (-N)v^{ii}(m,\,-m)+mv^{11}(N,\,-N) 
\end{equation*}
by direct computation, it follows that 
\begin{align*}
v^{1i}(-N,\,m)v^{1i}(N,\,-m)u & = 
  \bigl\{(-N)v^{ii}(m,\,-m)+mv^{11}(N,\,-N)\bigr\}u \\ 
  & \hspace*{40mm} + v^{1i}(N,\,-m)v^{1i}(-N,\,m)u. 
\end{align*}
Here we note that $\lambda+\Lambda_{i,m}-\Lambda_{1,N}$ 
is not contained in $Q_{+}$ since $\lambda(h_{1,N})=0$. 
Because the weight of $v^{1i}(-N,\,m)u$ is equal to 
$\lambda+\Lambda_{i,m}-\Lambda_{1,N} \notin Q_{+}$ 
by Remark~\ref{rem:wt}, we see from \eqref{eq:wsd} 
that $v^{1i}(-N,\,m)u=0$. Hence we get 
\begin{align*}
v^{1i}(-N,\,m)v^{1i}(N,\,-m)u & = 
  \bigl\{(-N)v^{ii}(m,\,-m)+mv^{11}(N,\,-N)\bigr\}u \\
  & = mN \underbrace{\lambda(h_{i,m})}_{\ge 1} u-
      mN \underbrace{\lambda(h_{1,N})}_{=0}u \ne 0.
\end{align*}
Thus we obtain $v^{1i}(-N,\,m)v^{1i}(N,\,-m)u \ne 0$, 
which implies that $v^{1i}(N,\,-m)u \ne 0$. 
It follows from Remark~\ref{rem:wt} that 
the weight $\lambda'$ of 
$v^{1i}(N,\,-m)u$ is equal to
$\lambda-\Lambda_{i,m}+\Lambda_{1,N}$. 
Since $\theta(\lambda')=\theta(\lambda)-1$, 
there exists $x' \in U(\CL_{r})$ such that  
$x' v^{1i}(N,\,-m)u \in \Mon$ by the inductive assumption. 
Thus we have proved Lemma~\ref{lem4:2}. 
\end{proof}
%
%
\subsection{Proof of the ``only if'' part of 
Proposition~\ref{prop:step2}.}
\label{subsec:s2-onlyif}

Assume that $W \subset \Mon$ is an $\Lon$-submodule of $\Mon$. 
Then we see that $W$ is stable under the action of $\CH$; 
indeed, we have $h_{1,\,l}W \subset W$ 
for all $l \in \BZ_{< 0}$ by assumption, and 
$h_{k,\,l}W=\{0\}$ for all $2 \le k \le d$ and 
$l \in \BZ_{< 0}$ (see \eqref{eq:hkl} and Remark~\ref{rem:wt2}). 
Thus, $W$ also admits 
the weight space decomposition as follows:
$W=\bigoplus_{\lambda \in Q_{+}} W^{\lambda}$ with 
$W^{\lambda}=W \cap (\Mam)^{\lambda}$.
Since $W \subset M_{r}^{(1)}=U(\Lon)\vm$, 
we deduce from Remark~\ref{rem:wt} 
(see also Remark~\ref{rem:wt2}) that 
$W^{\lambda} = \{0\}$ unless 
$\lambda \in Q_{+}^{(1)}:=
\sum_{l \in \BZ_{< 0}} \BZ_{\ge 0}\,\Lambda_{1,\,l}$. 
Hence we have 
%
%
\begin{equation} \label{eq:wtsd-1}
W=\bigoplus_{\lambda \in Q_{+}} W^{\lambda}
 =\bigoplus_{\lambda \in Q_{+}^{(1)}} W^{\lambda}.
\end{equation}

Now, let us prove that if $\Mam$ is an irreducible $\Lr$-module,
then $\Mon$ is an irreducible $\Lon$-module. 
It suffices to show the next lemma.
%
%
\begin{lem}\label{lem4:4}
Assume that $\Mon$ is a reducible $\Lon$-module, and 
let $W \subset \Mon$ be a proper $\Lon$-submodule of $\Mon$. 
Let $u$ be a nonzero homogeneous element of $W$, that is, 
$u \in W^{\lambda} \stm \{0\}$ for some $\lambda \in Q_{+}^{(1)}$. 
Then the $\Lr$-module $U(\CL_r)u \ (\subset \Mam)$ generated 
by the $u$ is a proper submodule of $\Mam$. Therefore 
the $\Lr$-module $\Mam$ is reducible. 
\end{lem}

\begin{proof}
Suppose that $U(\CL_r)u$ coincides with the whole of $\Mam$. 
Then there exists $x \in U(\CL_r)$ such that 
$xu \in \BC\vm \stm \{0\}$. Let 
\begin{equation*}
\CB_{1}:= \CB_{-} \cup 
 \bigl\{v^{1j}(m,\,n) \in \CB \mid 
   2 \le j \le d,\,
   m \in \BZ,\,n \in \BZ_{< 0} \bigr\},
\qquad 
\CB_{2}:=\CB \stm \CB_{1}.
\end{equation*}
Since an element in $\BC \subset \CL_{r}$ 
acts as a scalar multiple, 
we may assume, by the Poincar\'e-Birkhoff-Witt theorem, 
that the $x \in U(\CL_{r})$ above is of the form: 
$x=\sum_{1 \le t \le s} \alpha_{t} y_{t}z_{t}$, 
where $y_{t}$ (resp., $z_{t}$) is a product of elements 
in $\CB_{1}$ (resp., $\CB_{2}$) for each $1 \le t \le s$, 
and $\alpha_{t} \in \BC$ for each $1 \le t \le s$. 
Because $u$ is a homogeneous element, 
we see from Remark~\ref{rem:wt} that 
$y_{t}z_{t}u$ are also homogeneous elements 
for all $1 \le t \le s$. Since $xu \in \BC\vm \stm \{0\}$, and 
since $(M_{r})^{0}=(M_{r})_{0}=\BC\vm$, 
it follows that $y_{t}z_{t}u \in \BC\vm \stm \{0\}$ 
for some $1 \le t \le s$. 
Thus we may assume from the beginning 
that $x$ is of the form: $x=yz$, where 
$y$ (resp., $z$) is a product of 
elements in $\CB_{1}$ (resp., $\CB_{2}$).

Suppose that $y \ne 1$. Because $xu=yzu \in 
\BC\vm \stm \{0\} \subset (M_{r})^{0}$, we deduce 
from the definition of the set $\CB_{1}$ 
and Remark~\ref{rem:wt} that 
the weight of $zu$ is not contained in $Q_{+}$.
Hence, $zu=0$ by \eqref{eq:wsd}, 
which is a contradiction. Thus we get $y=1$. 
Write $z=z_{p}z_{p-1} \cdots z_{1}$ with 
$z_{q} \in \CB_{2}$ for $1 \le q \le p$. 
Suppose that there exists $1 \le q' \le p$ 
such that $z_{q'} \notin \CB^{(1)}=
 \bigl\{v^{11}(m,n) \mid 
 m,\,n \in \BZ \text{ with } m \le n \bigr\}$. 
Let $q:=\min \bigl\{1 \le q' \le p \mid 
 z_{q'} \notin \CB^{(1)}\bigr\}$. 
Then, $z_{q}$ is either of the following form: 
$v^{1j}(m,\,n)$ for some $2 \le j \le d$ and 
$m \in \BZ$, $n \in \BZ_{\ge 0}$, or 
$v^{ij}(m,\,n) \in \CB_{+}$ for some $2 \le i \le j \le d$.
Because $z_{1},\,z_{2},\,\dots,\,z_{q-1} \in \CB^{(1)}$, 
it is obvious that $z_{q-1} \cdots z_{1}u \in W \subset \Mon$. 
Thus, using Lemma~\ref{lem:zero} below, we see that 
$z_{q}z_{q-1} \cdots z_{1}u=0$, which is a contradiction.
Thus we conclude that $z_{1},\,z_{2},\,\dots,\,z_{p} 
\in \CB^{(1)}$. Hence, 
\begin{align*}
xu & = yzu = zu, \qquad \text{since $y=1$} \\
   & = z_{p}z_{p-1} \cdots z_{1}u \in W, \qquad
     \text{since $z_{1},\,z_{2},\,\dots,\,z_{p} \in \CB^{(1)}$}.
\end{align*}
Since $xu \in \BC\vm \stm \{0\}$, it follows that $\vm \in W$, 
which implies that $W=\Mon$. However, this is a contradiction, 
since $W$ is assumed to be a proper $\Lon$-submodule of $\Mon$. 
Thus we have proved Lemma~\ref{lem4:4}.
\end{proof}

Let us show the following lemma, 
which has been used in the proof of Lemma~\ref{lem4:4}. 
%
%
\begin{lem} \label{lem:zero}
{\rm (1)} Let $2 \le j \le d$, and $m \in \BZ$, $n \in \BZ_{\ge 0}$. 
Then, $v^{1j}(m,\,n)\Mon=\{0\}$. 

\noindent
{\rm (2)} Let $v^{ij}(m,\,n) \in \CB_{+}$ with $2 \le i \le j \le d$. 
Then, $v^{ij}(m,\,n)\Mon=\{0\}$. 
\end{lem}

\begin{proof}
(1) Let $\CS^{(1)}$ be the subset of $\CS$ consisting of 
all finite sequences of elements in $\CB^{(1)}_{-}:=\CB^{(1)} \cap \CB_{-}$ 
that is weakly decreasing with respect to the total ordering $\succ$. 
Then, $\BB^{(1)}:=\bigl\{ w(\bx) \mid \bx \in \CS^{(1)} \bigr\} 
\subset \BB$ is a linear basis of $\Mon$. Therefore 
it suffices to show that $v^{1j}(m,\,n)w(\bx)=0$ 
for all $\bx \in \CS^{(1)}$. 
Let $\bx=(x_{p} \succeq x_{p-1} \succeq \cdots \succeq x_{1}) 
\in \CS^{(1)}$ with $x_{q} \in \CB^{(1)}$ for $1 \le q \le p$. 
We show $v^{1j}(m,\,n)w(\bx)=0$ by 
induction on the length $p$ of the sequence $\bx$. 
If $p=0$, then the claim is obvious since $w(\bx)=\vm$. 
Assume that $p > 0$. Then, 
\begin{align*}
v^{1j}(m,\,n)w(\bx) & = 
 v^{1j}(m,\,n)x_{p}x_{p-1} \cdots x_{1}\vm \\
 & = [v^{1j}(m,\,n),\,x_{p}]x_{p-1} \cdots x_{1}\vm + 
     x_{p} \bigl\{v^{1j}(m,\,n)x_{p-1} \cdots x_{1}\vm\bigr\}.
\end{align*}
Since $v^{1j}(m,\,n)x_{p-1} \cdots x_{1}\vm=0$
by the inductive assumption, 
the second term of the right-hand side is equal to $0$. 
Assume that $x_{p}=v^{11}(s,t)$ with $s \le t < 0$. 
By simple computation, we see that 
\begin{align*}
[v^{1j}(m,\,n),\,x_{p}]_{r} & = 
[v^{1j}(m,\,n),\,v^{11}(s,t)]_{r} \\
& = \delta_{m+s,0} m v^{1j}(t,\,n) + 
    \delta_{m+t,0} m v^{1j}(s,\,n).
\end{align*}
Hence it follows from the inductive assumption that 
$[v^{1j}(m,\,n),\,x_{p}]x_{p-1} \cdots x_{1}\vm=0$.
Thus we get $v^{1j}(m,\,n)w(\bx)=0$, thereby completing 
the proof of part (1). 

\noindent
(2) Since $2 \le i \le j \le d$, it can be easily seen that 
$[v^{ij}(m,\,n),\,x]=0$ for all $x \in \CB^{(1)}$. Also, 
$v^{ij}(m,\,n)\vm=0$, since $v^{ij}(m,\,n) \in \CB_{+}$. 
The assertion of part (2) follows immediately from these facts. 
This completes the proof of the lemma.
\end{proof}

Proposition~\ref{prop:step2} follows from 
the results obtained in \S\S\ref{subsec:s2-if} 
and \ref{subsec:s2-onlyif}. At the end of this section, 
we show the following proposition, 
which is needed in \S\ref{subsec:quot}. 
%
%
\begin{prn} \label{prop:max2}
Assume that the $\CL_{r}$-module $\Mam$ is reducible, 
or equivalently, the $\Lon$-module $\Mon$ is reducible. 

\noindent 
{\rm (1)} If $W_{1}$ is the maximal 
proper $\Lon$-submodule of $\Mon$, 
then $U(\CL_{r})W_{1}$ is the maximal proper 
$\CL_{r}$-submodule of $\Mam$.

\noindent
{\rm (2)} If $W_{2}$ is the maximal proper
$\CL_{r}$-submodule of $\Mam$, 
then $W_{2} \cap \Mon$ is the maximal proper 
$\Lon$-submodule of $\Mon$. 
\end{prn}

\begin{proof}
We show that $U(\CL_{r})W_{1}=W_{2}$ (for part (1)), and 
$W_{1}=W_{2} \cap \Mon$ (for part (2)). 
Then we deduce from Lemma~\ref{lem4:4} that 
$U(\CL_{r})W_{1}$ is an $\CL_{r}$-submodule of $\Mam$ 
such that $U(\CL_{r})W_{1} \ne \Mam$. 
Hence we have $U(\CL_{r})W_{1} \subset W_{2}$ 
by the maximality of $W_{2}$. 
Also it follows from Lemma~\ref{lem4:2} 
(and the comment after it) that $W_{2} \cap \Mon$ 
is an $\Lon$-submodule of $\Mon$ such that 
$W_{2} \cap \Mon \ne \Mon$. 
Hence we have $W_{2} \cap \Mon \subset W_{1}$ 
by the maximality of $W_{1}$. 
Thus we obtain
\begin{equation*}
U(\CL_{r})W_{1} \cap \Mon 
\subset W_{2} \cap \Mon 
\subset W_{1}.
\end{equation*}
Because it is obvious that 
$W_{1} \subset U(\CL_{r})W_{1} \cap \Mon$, 
we get 
\begin{equation*}
U(\CL_{r})W_{1} \cap \Mon = W_{2} \cap \Mon = W_{1}, 
\end{equation*}
which shows part (2). 

Next, let us show $U(\CL_{r})W_{1}=W_{2}$ (i.e., part (1)).
Since $U(\CL_{r})W_{1} \subset W_{2}$ as shown above, 
it suffices to show that $U(\CL_{r})W_{1} \supset W_{2}$. 
Note that $W_{2}$ admits the weight space decomposition 
$W_{2}=\bigoplus_{\lambda \in Q_{+}} (W_{2})^{\lambda}$ 
with $(W_{2})^{\lambda}=W_{2} \cap (\Mam)^{\lambda}$. 
Let $u \in W_{2}$ be a homogeneous element of 
weight $\lambda \in Q_{+}$, that is, 
$u \in (W_{2})^{\lambda}$. 
We show by induction on $\theta(\lambda)$ that 
$u \in U(\CL_{r})W_{1}$ 
(for the definition of $\theta(\lambda)$, 
see \eqref{eq:theta}). 
If $\theta(\lambda) = 0$, then $u \in \Mon$ 
by Remark~\ref{rem:wt2}, and hence $u \in W_{2} \cap \Mon$. 
Since $W_{2} \cap \Mon=W_{1}$ as shown above, it follows that
$u \in W_{1} \subset U(\CL_{r})W_{1}$. 
Next, let us assume that $\theta(\lambda) > 0$. 
Let $2 \le i \le d$ and $m \in \BZ_{< 0}$ 
be such that $\lambda(h_{i,m}) \ge 1$, and 
$N \in \BZ_{< 0}$ such that $\lambda(h_{1,N})=0$.
Then we deduce from the proof of 
Lemma~\ref{lem4:2} that 
$v^{1i}(N,\,-m)u \ne 0$, and the weight $\lambda'$ of 
$v^{1i}(N,\,-m)u$ satisfies $\theta(\lambda')=\theta(\lambda)-1$. 
Hence it follows from the inductive assumption that 
$v^{1i}(N,\,-m)u \in U(\CL_{r})W_{1}$. 
Further, as in the proof of Lemma~\ref{lem4:2}, 
we deduce that 
\begin{equation*}
v^{1i}(-N,\,m) 
\underbrace{v^{1i}(N,\,-m)u}_{\in U(\CL_{r})W_{1}} =  
mN \underbrace{\lambda(h_{i,\,m})}_{\ge 1} u,
\end{equation*}
which implies that $u \in U(\CL_{r})W_{1}$. 
This completes the proof of Proposition~\ref{prop:max2}.
\end{proof}
%
%
\section{Irreducibility of $\Mon$
for $r \in \BC \stm \BZ$.}
\label{sec:step3}

This subsection is devoted to proving 
the following proposition.
%
%
\begin{prn} \label{prop:step3}
If $r \in \BC \stm \BZ$, 
then $\Mon$ is an irreducible $\Lon$-module.
\end{prn}

%
\subsection{Notation and some lemmas.}
\label{subsec:notation}
For simplicity of notation, we set $v(m,\,n):=v^{11}(m,\,n)$ 
for $m,\,n \in \BZ$, and 
$\Lambda_{l}:=\Lambda_{1,l} \in \CH^{\ast}$ 
for $l \in \BZ_{< 0}$. 
Recall that $\CS^{(1)}$ denotes the subset of $\CS$ consisting of 
all finite sequences of elements in 
$\CB^{(1)}_{-}=\CB^{(1)} \cap \CB_{-}$ that 
is weakly decreasing with respect to the total ordering $\succ$, 
and $\BB^{(1)}=
 \bigl\{ w(\bx) \mid \bx \in \CS^{(1)} \bigr\} \subset \BB$ 
is a linear basis of $\Mon$. For each $\lambda \in Q_{+}^{(1)}$, 
we denote by $\BB^{(1)}_{\lambda}$ 
the set of all elements in $\BB^{(1)}$ 
whose weight is equal to $\lambda$, and set 
$\CS^{(1)}_{\lambda}:=\bigl\{\bx \in \CS^{(1)} \mid 
w(\bx) \in \BB^{(1)}_{\lambda} \bigr\}$. 

Let $u \in \Mon$, and write it as a linear combination of 
elements of $\BB^{(1)}$: 
$u=\sum_{b \in \BB^{(1)}} \alpha_{b}\hspace{1pt}b$ 
with $\alpha_{b} \in \BC$ for $b \in \BB^{(1)}$. 
Then we set $\BB[u]:=
\bigl\{b \in \BB^{(1)} \mid \alpha_{b} \ne 0\bigr\}$, and 
$\CS[u]:=\bigl\{\bx \in \CS^{(1)} \mid 
w(\bx) \in \BB[u]\bigr\}$. 

The following formulas can be 
shown by simple computation. 
%
%
\begin{lem} \label{lem:LB}
{\rm (1)} 
Let $s,\,t \in \BZ_{> 0}$, and $m,\,n \in \BZ_{> 0}$. 
Then, 
%
%
\begin{equation} \label{eq:LB02}
[v(-m,\,n),\,v(-s,\,-t)]_{r} = 
n \bigl\{\delta_{n,s}v(-m,\,-t)+\delta_{n,t}v(-s,\,-m)\bigr\}.
\end{equation}
{\rm (2)} 
Let $s,\,t \in \BZ_{> 0}$ with $s \ne t$, and 
$m \in \BZ_{> 0}$. Then, 
%
%
\begin{equation} \label{eq:LB03}
[v(m,\,m),\,v(-s,\,-t)]_{r} = 
2m \bigl\{\delta_{m,s}v(-t,\,m)+\delta_{m,t}v(-s,\,m)\bigr\}.
\end{equation}
\end{lem}

%
\begin{lem} \label{lem:vmm}
Let $m \in \BZ_{> 0}$, and $\nu \in \BZ_{\ge 0}$. Then, 
%
%
\begin{equation} \label{eq:vmm}
v(m,\,m)v(-m,\,-m)^{\nu}\vm
 =2m^{2}\nu(r+2\nu-2)v(-m,\,-m)^{\nu-1}\vm.
\end{equation}
\end{lem}
%
%
\subsection{Proof of Proposition~\ref{prop:step3}.}
\label{subsec:prf-s3}
We show that if $W \subset \Mon$ is 
a nonzero $\Lon$-submodule of $\Mon$, then $W=\Mon$. 
Suppose that $W \subsetneq \Mon$.
Because $W$ admits the weight space decomposition 
$W=\bigoplus_{\lambda \in Q_{+}^{(1)}} W^{\lambda}$ 
with respect to $\CH$ (see \S\ref{subsec:s2-onlyif}), 
it can be easily seen by usual way 
(see also Remark~\ref{rem:wt}) that 
$W$ contains a (homogeneous) singular vector $u$, i.e., 
a nonzero element $u$ such that $u \in W^{\lambda}$ 
for some $\lambda \in Q_{+}^{(1)} \stm \{0\}$, and 
$v(m,\,n)u=0$ for all $m,\,n \in \BZ$ with 
$m+n > 0$. 
%
%
\begin{claim} \label{c:s3-1}
The set $\bigl\{ l \in \BZ_{< 0} \mid \lambda(h_{l}) > 0 \bigr\}$
is identical to $\bigl\{-p,\,-p+1,\,\dots,\,-2,\,-1\bigr\}$ 
for some $p \in \BZ_{\ge 1}$. 
\end{claim}

\noindent 
{\it Proof of Claim~\ref{c:s3-1}.} 
Assume that $\bigl\{ l \in \BZ_{< 0} \mid \lambda(h_{l}) > 0 \bigr\} = 
\bigl\{l_{p} < l_{p-1} < \cdots < l_{1}\bigr\}$; note that 
$l_{q} \le -q$ for every $1 \le q \le p$. Suppose that 
$l_{q} < -q$ for some $1 \le q \le p$. We set 
$q_{0}:=\min \bigl\{ 1 \le q \le p \mid l_{q} < -q \bigr\}$. 
Since $\lambda(h_{l_{q_0}}) > 0$ and $\lambda(h_{-q_0})=0$, 
we deduce by a way similar to the proof of 
Lemma~\ref{lem4:2} that $v(l_{q_{0}},\,q_{0})u=0$, and 
\begin{align*}
v(l_{q_{0}},\,q_{0})v(-q_{0},\,-l_{q_{0}})u 
& = l_{q_0}v(-q_{0},q_{0})u+q_{0}v(l_{q_{0}},\,-l_{q_{0}})u \\
& = q_{0}l_{q_0} \underbrace{\lambda(h_{-q_0})}_{=0} u+ 
    (-q_{0}l_{q_0}) \underbrace{\lambda(h_{l_{q_0}})}_{> 0} u \ne 0. 
\end{align*}
Hence we get $v(-q_{0},\,-l_{q_{0}})u \ne 0$. 
However, since $l_{q_{0}} < -q_{0}$, this contradicts the assumption 
that $u$ is a singular vector. Thus we obtain 
$l_{q}=-q$ for all $1 \le q \le p$, thereby completing 
the proof of Claim~\ref{c:s3-1}. 
%
%
\begin{claim} \label{c:s3-2}
The weight $\lambda$ of the singular vector $u$ 
is of the form\,{\rm:} 
$\lambda=\sum_{q=1}^{p} 2\nu_{q}\Lambda_{-q}$ 
with $\nu_{q} \in \BZ_{> 0}$ for $1 \le q \le p$. 
Furthermore, the set $\CS[u]$ contains the element 
$\bx_{\lambda} \in \CS^{(1)}_{\lambda}$ such that 
\begin{equation*}
w(\bx_{\lambda})= 
 \prod_{1 \le q \le p} v(-q,\,-q)^{\nu_{q}} \vm \in \BB^{(1)}_{\lambda}. 
\end{equation*}
\end{claim}

\noindent 
{\it Proof of Claim~\ref{c:s3-2}.} 
For $\bx \in \CS[u]$ and $1 \le q \le p$, 
we define $\kappa_{q}(\bx)$ to be the number of 
$v(-q,\,-q)$ appearing in the sequence $\bx$. 
Take $\bx \in \CS[u]$ such that the sum 
$\sum_{q=1}^{p} \kappa_{q}(\bx)$ is maximum, and 
assume that $w(\bx) \in \BB[u]$ is of the form:
\begin{equation*}
w(\bx) = 
  \prod_{1 \le q \le p}
  v(-q,\,-q)^{\nu_{q}}
  \prod_{1 \le t < s \le p} v(-s,\,-t)^{\nu_{s,t}} \vm
\end{equation*}
for some $\nu_{q} \in \BZ_{\ge 0}$, $1 \le q \le p$, 
and $\nu_{s,t} \in \BZ_{\ge 0}$, $1 \le t < s \le p$. 
Suppose that $\nu_{s_{1},t_{1}} > 0$ for some 
$1 \le t_{1} < s_{1} \le p$.
Then we prove that $v(-t_{1},\,s_{1})u \ne 0$.
For this, it suffices to show that the set 
$\BB[v(-t_{1},\,s_{1})u]$ is not empty. 
Here we show that $\BB[v(-t_{1},\,s_{1})u]$ 
contains the following element: 
\begin{align*}
\prod_{1 \le q \le p,\,q \ne t_{1}}
  v(-q,\,-q)^{\nu_{q}} \, 
  v(-t_{1},\,-t_{1})^{\nu_{t_{1}}+1} 
\prod_{
  \begin{subarray}{c}
  1 \le t < s \le p, \\[1mm]
  (s,t) \ne (s_{1},t_{1})
  \end{subarray}
  } v(-s,\,-t)^{\nu_{s,t}}\,
  v(-s_{1},\,-t_{1})^{\nu_{s_{1},t_{1}}-1}\,\vm.
\end{align*}
It can be easily seen by \eqref{eq:LB02} that 
the element above is contained in 
$\BB[v(-t_{1},\,s_{1})w(\bx)]$. On the other hand, 
we deduce, using \eqref{eq:LB02} and 
the maximality of $\bx$, that the element above is not 
contained in $\BB[v(-t_{1},\,s_{1})w(\by)]$ for all 
$\by \in \CS[u]$ with $\by \ne \bx$. Combining these, 
we conclude that the element above is contained in 
$\BB[v(-t_{1},\,s_{1})u]$, and hence 
$v(-t_{1},\,s_{1})u \ne 0$. However, since $t_{1} < s_{1}$, 
this contradicts the assumption that 
$u$ is a singular vector. 
Thus we obtain $\nu_{s,t}=0$ for all $1 \le t < s \le p$, 
thereby completing the proof of Claim~\ref{c:s3-2}. 

\vspace{3mm}

By Claim~\ref{c:s3-2}, we may assume that 
the singular vector $u \in W^{\lambda}$ 
is of the form: 
\begin{equation*}
u = w(\bx_{\lambda}) + 
    \sum_{1 \le t < s \le p}
    \alpha_{s,t} w(\bx_{\lambda}^{s,t}) + 
    \sum_{\bx \in \CS'} \beta_{\bx} w(\bx)
\end{equation*}
with $\alpha_{s,t} \in \BC$ for $1 \le t < s \le p$, and 
$\beta_{\bx} \in \BC$ for $\bx \in \CS'$, where 
$\bx_{\lambda}^{s,t}$ is 
the element of $\CS^{(1)}_{\lambda}$ such that 
\begin{equation*}
w(\bx_{\lambda}^{s,t})= 
  v(-s,\,-s)^{\nu_{s}-1}\,
  v(-t,\,-t)^{\nu_{t}-1}\,
  v(-s,\,-t)^{2}
 \prod_{1 \le q \le p,\,q \ne t,\,s} 
 v(-q,\,-q)^{\nu_{q}} \vm \in \BB^{(1)}_{\lambda}, 
\end{equation*}
and $\CS':=\CS^{(1)}_{\lambda} \stm 
\{\bx_{\lambda},\,\bx_{\lambda}^{s,t} \mid 1 \le t < s \le p\bigr\}$. 
%
%
\begin{claim} \label{c:s3-3}
{\rm (1)} We have 
$\alpha_{s,t}=-\nu_{s}$ for every $1 \le t < s \le p$.

\noindent
{\rm (2)} For every $1 \le s \le p$, 
\begin{equation*}
\nu_{s} (r+2\nu_{s}-2) + 
 \sum_{1 \le t \le s-1} \alpha_{s,t} + 
 \sum_{s+1 \le t \le p} \alpha_{t,s}=0.
\end{equation*}
\end{claim}

\noindent 
{\it Proof of Claim~\ref{c:s3-3}.} 
(1) Fix $1 \le t < s \le p$. 
We can easily check by using \eqref{eq:LB02} that 
if the set $\BB[v(-t,\,s)w(\bx)]$ with 
$\bx \in \CS^{(1)}_{\lambda}$ contains 
\begin{equation*}
w_{1}:=
 v(-s,\,-s)^{\nu_{s}-1}\,
 v(-s,\,-t)
 \prod_{1 \le q \le p,\,q \ne s} 
 v(-q,\,-q)^{\nu_{q}} \vm, 
\end{equation*}
then $\bx=\bx_{\lambda}$ or 
$\bx=\bx_{\lambda}^{s,t}$. 
By simple computation, along with \eqref{eq:LB02}, 
we get 
\begin{align*}
& v(-t,\,s) w(\bx_{\lambda}) = 2\nu_{s}s w_{1}, & 
& v(-t,\,s) w(\bx_{\lambda}^{s,t}) = 2s w_{1} + \text{(other term)}.
\end{align*}
Because $v(-t,\,s)w=0$ by assumption, it follows that 
$2\nu_{s}s+2s\alpha_{s,t}=0$, and hence $\alpha_{s,t}=-\nu_{s}$. 

\noindent
(2) Fix $1 \le s \le p$. We can easily check by using 
\eqref{eq:LB01}, \eqref{eq:LB02}, and \eqref{eq:LB03} that 
if the set $\BB[v(s,\,s)w(\bx)]$ with 
$\bx \in \CS^{(1)}_{\lambda}$ contains 
\begin{equation*}
w_{2}:=
 v(-s,\,-s)^{\nu_{s}-1}
 \prod_{1 \le q \le p,\,q \ne s} 
 v(-q,\,-q)^{\nu_{q}} \vm,
\end{equation*}
then $\bx=\bx_{\lambda}$, or 
$\bx=\bx_{\lambda}^{s,t}$ for $1 \le t \le s-1$, or 
$\bx=\bx_{\lambda}^{t,s}$ for $s+1 \le t \le p$. 
We see from Lemma~\ref{lem:vmm} that 
\begin{equation*}
v(s,\,s) w(\bx_{\lambda}) = 
 2s^{2}\nu_{s}(r +2\nu_{s}-2) w_{2}. 
\end{equation*}
Also it follows from \eqref{eq:LB03} that 
for each $1 \le t \le s-1$, 
%
%
\begin{align}
v(s,\,s)  w(\bx_{\lambda}^{s,t}) 
  & = 2s X_{1} v(s,\,-t)\,v(-s,\,-t)\,v(-s,\,-s)^{\nu_{s}-1} \vm \notag \\
  & \qquad +
      2s X_{1} v(-s,\,-t)\,v(s,\,-t)\,v(-s,\,-s)^{\nu_{s}-1} \vm \notag \\
  & \qquad +
      2s X_{1} v(-s,\,-t)^{2}v(s,\,s)v(-s,\,-s)^{\nu_{p}-1} \vm. \label{eq:3-4-1}
\end{align}
Here, for simplicity of notation, we set 
\begin{equation*}
X_{1}:=
 \prod_{1 \le q \le p,\,q \ne t,\,s} v(-q,\,-q)^{\nu_{q}}\,
 v(-t,\,-t)^{\nu_{s}-1}. 
\end{equation*}
The second and third terms of the right-hand side of 
\eqref{eq:3-4-1} do not contribute the coefficient of 
$w_{2}$ since they contain $v(-s,\,-t)$. 
Also, the first term is: 
\begin{align*}
& 2s X_{1} v(s,\,-t)\,v(-s,\,-t)\,v(-s,\,-s)^{\nu_{s}-1} \vm \\
& \hspace*{30mm} = 
  2s^{2} X_{1} v(-t,\,-t)\,v(-s,\,-s)^{\nu_{p}-1} \vm + 
  \text{(other term)}.
\end{align*}
Thus we obtain 
\begin{equation*}
v(s,\,s) w(\bx_{\lambda}^{s,t}) = 
  2s^{2} w_{2} + \text{(other terms)}
\end{equation*} 
for $1 \le t \le s-1$. 
Similarly, we can show that 
\begin{equation*}
v(s,\,s) w(\bx_{\lambda}^{t,s}) = 
  2s^{2} w_{2} + \text{(other terms)}
\end{equation*}
for $s+1 \le t \le p$. 
Because $v(p,\,p)u=0$ by assumption, we obtain
\begin{equation*}
2s^{2}\nu_{s}(r +2\nu_{s}-2) + 
 2s^{2} \sum_{1 \le t \le s-1} \alpha_{s,t}+
 2s^{2} \sum_{s+1 \le t \le p} \alpha_{t,s} = 0,
\end{equation*}
and hence the equation of part (2). 
Thus we have proved Claim~\ref{c:s3-3}. 

\vspace{3mm}

Combining the equations in Claim~\ref{c:s3-3} 
with $s=p$, we obtain 
$\nu_{p}(r+2\nu_{p}-2) -(p-1)\nu_{p}=0$,
and hence that 
%
%
\begin{equation} \label{eq:rnup}
r+2\nu_{p}-1-p=0.
\end{equation}
However, this is a contradiction, 
since $r$ is assumed not to be an integer. 
This completes the proof of Proposition~\ref{prop:step3}. \qed
%
%
\begin{rmk} \label{rem:nu}
For later use, let us show the following assertion: 
Keep the notation in the proof of 
Proposition~\ref{prop:step3} above. Then, 
$\nu_{1}=\nu_{2}= \cdots = \nu_{p}$. 
Indeed, by the equations of Claim~\ref{c:s3-3}, 
we see that 
$\nu_{s}(r+2\nu_{s}-2)
 - (s-1)\nu_{s}-\nu_{s+1} - \cdots - \nu_{p}=0$ 
for every $1 \le s \le p$. Therefore, 
\begin{align*}
& (\nu_{s}-\nu_{s+1})(r+2\nu_{s+1}-s-1+2\nu_{s}) \\
& \qquad = 
  \Bigl\{
    \nu_{s}(r+2\nu_{s}-2)
    - (s-1)\nu_{s}-\nu_{s+1} - \cdots - \nu_{p}
  \Bigr\} \\
& \hspace*{20mm} - 
  \Bigl\{
    \nu_{s+1}(r+2\nu_{s+1}-2)
    - s\nu_{s+1}-\nu_{s+2} - \cdots - \nu_{p}
  \Bigr\} \\
& \qquad = 0
\end{align*}
for every $1 \le s \le p-1$. Using this equation and 
\eqref{eq:rnup}, we can show by descending induction that 
$\nu_{q}=\nu_{q+1}$ for all $1 \le q \le p-1$, and hence 
$\nu_{1}=\nu_{2}= \cdots = \nu_{p}$.
\end{rmk}
%
%
\section{Reducibility of $\Mon$ for $r \in \BZ$.}
\label{sec:step4}
%
%
%
\subsection{Notation and proposition.}
\label{subsec:prop}
For each $p \in \BZ_{\ge 0}$, let $\bV_{p}$ be 
the following matrix of size $p$ with entries in $\CB^{(1)}_{-}$: 
\begin{equation*}
\bV_{p}=\bigl( v(-s,-t) \bigr)_{1 \le s,\,t \le p}.
\end{equation*}
Since $xy=yx$ for all $x,\,y \in \CB^{(1)}_{-}$, 
we can consider the determinant $\det \bV_{p}$ of 
the matrix $\bV_{p}$;
\begin{equation*}
\det \bV_{p}=\sum_{\sigma \in \FS_{p}} 
 \sgn (\sm) \prod_{1 \le q \le p}
 v(-q,\,-\sigma(q)),
\end{equation*}
where $\FS_{p}$ denotes the symmetric group of 
degree $p$, and $\sgn (\sm)$ denotes the signature of 
a permutation $\sm \in \FS_{p}$. 
In this subsection, we show 
the following proposition.
%
%
\begin{prn} \label{prop:step4}
Assume that $r \in \BZ$. 
Let $\nu \in \BZ_{\ge 1}$ and $p \in \BZ_{\ge 1}$
be positive integers satisfying the relation 
$r=1-2\nu+p$. Then, $(\det \bV_{p})^{\nu}\vm$ is 
a singular vector of $\Mon$, that is, 
$v(m,n)(\det \bV_{p})^{\nu}\vm=0$ 
for all $m,\,n \in \BZ$ with $m+n > 0$. 
Therefore the $\Lon$-module $\Mon$ is reducible. 
\end{prn}

Theorem~\ref{thm:main} 
follows immediately from Propositions~\ref{V=M}, 
\ref{prop:step1}, \ref{prop:step2}, \ref{prop:step3}, 
and \ref{prop:step4} (see also the comment 
after Theorem~\ref{thm:main}). 

%
\subsection{Proof of Proposition~\ref{prop:step4}.}
\label{subsec:proof-prop}
Let us first show the following lemmas.
%
%
\begin{lem} \label{lem:det01}
Let $m \in \BZ$ with $1 \le m \le p$, and 
$n \in \BZ_{\ge 0}$ with $m \ne n$. Then,
%
%
\begin{equation} \label{eq:det01}
[v(-m,\,n),\,\det \bV_{p}]=0.
\end{equation}
\end{lem}

\begin{proof}
If $n \ge p+1$ or $n=0$, then the assertion is obvious, since 
$[v(-m,\,n),\, v(-s,-t)]_{r}=0$ for all $1 \le s,\,t \le p$. 
Assume that $1 \le n \le p$. It can be seen from \eqref{eq:LB02} that 
$[v(-m,\,n),\, v(-s,-t)]_{r}$ is equal to $0$, 
or is contained in $\CB^{(1)}_{-}$, 
up to a scalar multiple. Since $\ad x$ is 
a derivation on $U(\Lon)$ for $x \in \Lon$, 
we deduce that 
\begin{equation*}
[v(-m,\,n),\,\det \bV_{p}]=
 \det \bW_{1} + \det \bW_{2} + \cdots + \det \bW_{p}, 
\end{equation*}
where $\bW_{s}$, $1 \le s \le p$, 
is the matrix obtained by replacing 
the $s$-th row $(v(-s,-t))_{1 \le t \le p}$ of the matrix 
$\bV_{p}$ with 
\begin{equation*}
\Bigl([v(-m,\,n),\, v(-s,-t)]_{r}\Bigr)_{1 \le t \le p}.
\end{equation*}
It follows from \eqref{eq:LB02} that 
if $s \ne n$, then 
\begin{equation*}
\Bigl([v(-m,\,n),\, v(-s,-t)]_{r}\Bigr)_{1 \le t \le p} 
 = \bigl(0,\,\dots,\,0,\,nv(-s,-m),\,0,\,\dots,\,0\bigr),
\end{equation*}
where $nv(-s,-m)$ is placed at the $n$-th entry. Also, 
it follows from \eqref{eq:LB02} that 
\begin{align*}
& \Bigl([v(-m,\,n),\, v(-n,-t)]_{r}\Bigr)_{1 \le t \le p} 
  - n \underbrace{\bigl(v(-m,-t)\bigr)_{1 \le t \le p}}_{%
      \text{the $m$-th row of $\bW_{n}$}} \\[3mm]
& \hspace*{45mm}
  = \bigl(0,\,\dots,\,0,\,nv(-n,-m),\,0,\,\dots,\,0\bigr),
\end{align*}
where $nv(-n,-m)$ is placed at the $n$-th entry. 
Thus, 
$[v(-m,\,n),\,\det \bV_{p}] = n \det \bV_{p}'$,
where $\bV_{p}'$ is the matrix obtained by replacing 
the $n$-th column $(v(-s,-n))_{1 \le s \le p}$ of 
the matrix $\bV_{p}$ with $(v(-s,-m))_{1 \le s \le p}$ 
(i.e., the $m$-th column  of $\bV_{p}$). 
Because the $n$-th and $m$-th columns of $\bV_{p}'$ are equal, 
we get $\det \bV_{p}'=0$, thereby completing 
the proof of the lemma.
\end{proof}
%
%
\begin{lem} \label{lem:det02}
Let $1 \le m \le p$, and assume that 
$u \in \Mon$ satisfies the conditions that 
$v(-t,\,m)u=0$ for all $1 \le t \le p$ with $t \ne m$, and 
$v(-m,\,m)u=\alpha m u$ for some $\alpha \in \BC$. 
Then we have 
\begin{equation*}
v(m,\,m) (\det \bV_{p})\,u 
    = 2m^{2}(2\alpha+r-p+1) \det (\bV_{p}^{(m)})u + 
      (\det \bV_{p})v(m,\,m)u,
\end{equation*}
where $\bV_{p}^{(m)}$ is the matrix obtained by 
removing the $m$-th row and the $m$-th column from 
the matrix $\bV_{p}$. 
\end{lem}

\begin{proof}
We set $\FS_{p-1}^{(m)}:=
 \bigl\{ \sigma \in \FS_{p} \mid \sigma(m)=m \bigr\}$, 
which is a subgroup of $\FS_{p}$ isomorphic to 
the symmetric group $\FS_{p-1}$ of degree $p-1$. 
Let 
\begin{equation*}
R^{(m)}:=\bigl\{1,\,(m,t) 
 \mid \text{$1 \le t \le p$ with $t \ne m$}\bigr\} 
\end{equation*}
be a complete set of the representatives for 
the quotient $\FS_{p}/\FS_{p-1}^{(m)}$, 
where $(m,t) \in \FS_{p}$ denotes the transposition 
interchanging $m$ and $t$. 
Decompose $\FS_{p}$ as: 
$\FS_{p}=\bigsqcup_{\tau \in \FS_{p-1}^{(m)}} \tau R^{(m)}$, 
with $\tau R^{(m)} := \bigl\{ \tau,\,\tau \cdot (m,t) \mid 
\text{$1 \le t \le p$ with $t \ne m$}\bigr\}$. Then we have
%
%
\begin{align} 
& v(m,\,m) (\det \bV_{p})\,u 
    = [v(m,\,m),\,\det \bV_{p}]u + 
      (\det \bV_{p})v(m,\,m)u \notag \\[3mm]
& \qquad = 
 \sum_{\tau \in \FS_{p-1}^{(m)}}
 \sum_{\sigma \in \tau R^{m}} \sgn(\sigma) 
 \left[v(m,\,m),\ \prod_{1 \le s \le p} v(-s,\,-\sigma(s))\right]u 
 + (\det \bV_{p})v(m,\,m)u. 
 \label{eq:det02-2}
\end{align}
If $\tau \in \FS_{p-1}^{(m)}$, then we see, 
using \eqref{eq:LB01}, \eqref{eq:LB03}, 
and the assumption on $u$, that 
\begin{align*}
& \left[v(m,\,m),\ \prod_{1 \le s \le p} v(-s,\,-\tau(s))\right]u  \\[3mm]
& \hspace*{20mm} =
  \left\{
  \prod_{1 \le s \le p,\,s \ne m}
  v(-s,\,-\tau(s))
  \right\} 
  [v(m,\,m),\, v(-m,\,-m)] u \\[3mm]
& \hspace*{20mm}
  = \left\{
    \prod_{1 \le s \le p,\,s \ne m} 
    v(-s,\,-\tau(s))
    \right\}
    \bigl\{ 4mv(-m,\,m)+2rm^{2} \bigr\}u \\[3mm]
& \hspace*{20mm}
  = 2m^{2}(2\alpha+r)
    \prod_{1 \le s \le p,\,s \ne m} 
    v(-s,\,-\tau(s))u.
\end{align*}
Assume that $\sigma=\tau \cdot (m,\,t) \in \tau R^{(m)}$ 
with $\tau \in \FS_{p-1}^{(m)}$ and $1 \le t \le p$ with $t \ne m$; 
note that $\sigma^{-1}(m)=t$, 
and $\sigma(s)=\tau(s)$ for $1 \le s \le p$ with $s \ne m,\,t$. 
Then, 
\begin{align*}
& \left[v(m,\,m),\ \prod_{1 \le s \le p} v(-s,\,-\sigma(s))\right]u \\
& \quad 
  = \left\{
     \prod_{\begin{subarray}{c}
     1 \le s \le p, \\[1mm]
     s \ne m,\,\sigma^{-1}(m)
     \end{subarray}} 
     v(-s,\,-\sigma(s))
   \right\} 
   \Bigl[v(m,\,m),\ v(-\sigma^{-1}(m),\,-m)v(-m,\,-\sigma(m))\Bigr] u \\[3mm]
& \quad 
  = \left\{
     \prod_{
     1 \le s \le p,\,
     s \ne m,\,t} 
     v(-s,\,-\tau(s))
   \right\}
   \Bigl[v(m,\,m),\ v(-t,\,-m) v(-m,\,-\tau(t))\Bigr]u. 
\end{align*}
Using \eqref{eq:LB02}, \eqref{eq:LB03}, 
and the assumption for $u$ 
(note that $\tau(t) \ne m$ and $t \ne m$), 
we deduce that 
\begin{equation*}
\Bigl[v(m,\,m),\ v(-t,\,-m)v(-m,\,-\tau(t))\Bigr]u =
 2m^{2}v(-t,\,-\tau(t))u.
\end{equation*}
Thus, 
\begin{equation*}
\left[v(m,\,m),\ \prod_{1 \le s \le p} v(-s,\,-\sigma(s))\right]u = 
2m^{2} \prod_{1 \le s \le p,\,s \ne m} v(-s,\,-\tau(s))u.
\end{equation*}
Therefore, for each $\tau \in \FS_{p-1}^{(m)}$, 
\begin{align*}
 & \sum_{\sigma \in \tau R^{(m)}} \sgn(\sigma) 
   \left[v(m,\,m),\ \prod_{1 \le s \le p} v(-s,\,-\sigma(s))\right]u  \\[3mm]
 & = 2m^{2}
     \left\{
       (2\alpha+r) \sgn (\tau) + 
       \sum_{1 \le t \le p,\,t \ne m}
       \sgn(\tau \cdot (1,t))
     \right\}
     \prod_{1 \le s \le p,\, s \ne m} v(-s,\,-\tau(s))u \\[3mm]
 & = 2m^{2}(2\alpha+r-p+1) \, 
     \sgn (\tau)
     \prod_{1 \le s \le p,\, s \ne m} v(-s,\,-\tau(s))u.
\end{align*}
Combining this equation and equation \eqref{eq:det02-2}, 
we see that 
\begin{equation*}
v(m,\,m) (\det \bV_{p})\,u
    = 2m^{2}(2\alpha+r-p+1) \det (\bV_{p}^{(m)})u + 
      (\det \bV_{p})v(m,\,m)u,
\end{equation*}
as desired. 
\end{proof}

\begin{proof}[Proof of Proposition~\ref{prop:step4}]
It is obvious that $(\det \bV_{p})^{\nu}\vm$ is 
a nonzero homogeneous element, and 
$(\det \bV_{p})^{\nu}\vm \notin \BC\vm$. 
Let $m,\,n \in \BZ$ be such that $m \le n$ and $m+n > 0$; 
note that $n \ge 1$. 
If $n \ge p+1$, then $v(m,\,n) (\det \bV_{p})^{\nu}\vm=0$ 
since the weight of $v(m,\,n) (\det \bV_{p})^{\nu}\vm=0$ is 
not contained in $Q_{+}^{(1)}$ (see \eqref{eq:wtsd-1}). 
Let us consider the case that $1 \le n \le p$. 

\paragraph{Case 1.} 
Assume that $1 \le n \le p$ and $-n < m < 0$. 
It follows from Lemma~\ref{lem:det01} that 
$[v(m,\,n),\,\det \bV_{p}]=0$. Hence
\begin{equation*}
v(m,\,n) (\det \bV_{p})^{\nu}\vm = 
(\det \bV_{p})^{\nu}v(m,\,n)\vm=0.
\end{equation*}

\paragraph{Case 2.} 
Assume that $1 \le n \le p$ and $m=0$. 
By direct computation, we see that 
$v(0,\,n)=(p+1)^{-1}[v(n,\,p+1),\,v(-p-1,\,0)]_{r}$. 
As mentioned above, we have 
$v(n,\,p+1)(\det \bV_{p})^{\nu}\vm=0$. 
Also, since $[v(-s,-t),\,v(-p-1,\,0)]_{r}=0$ 
for all $1 \le s,\,t \le p$, it follows that 
$[v(-p-1,\,0),\,\det \bV_{p}]=0$, and hence that 
$v(-p-1,\,0) (\det \bV_{p})^{\nu}\vm = 0$ as above.
Thus we get $v(0,\,n)(\det \bV_{p})^{\nu}\vm=0$. 

\paragraph{Case 3.}
Assume that $1 \le n \le p$ and $1 \le m \le n$. 
First let us consider the case that $m=n$. 
Namely, we show that 
$v(m,\,m)(\det \bV_{p})^{\nu}\vm=0$ 
for all $1 \le m \le p$. 
It follows immediately 
from Lemma~\ref{lem:det01} that 
$u=(\det \bV_{p})^{\nu_{1}}\vm$, $\nu_{1} \in \BZ_{\ge 0}$, 
satisfies the assumption of Lemma~\ref{lem:det02} 
with $\alpha=2\nu_{1}$. 
Hence, by using Lemma~\ref{lem:det02} repeatedly, 
we obtain
\begin{align*}
& v(m,\,m) (\det \bV_{p})^{\nu}\vm = 
  v(m,\,m) (\det \bV_{p}) (\det \bV_{p})^{\nu-1}\vm \\
& = 
  2m^{2}\{4(\nu-1)+r-p+1\}(\det \bV_{p}^{(m)})(\det \bV_{p})^{\nu-1}\vm + 
  (\det \bV_{p})v(m,\,m)(\det \bV_{p})^{\nu-1}\vm \\
& = \cdots \cdots \\
& = 2m^{2}(\det \bV_{p}^{(m)})\left\{
    \sum_{\nu_{1}=0}^{\nu-1}\bigl(4\nu_{1}+r-p+1\bigr)\right\}
    (\det \bV_{p})^{\nu-1}\vm \\[3mm]
& = 4m^{2}\nu(2\nu-1+r-p)(\det \bV_{p}^{(m)})(\det \bV_{p})^{\nu-1}\vm. 
\end{align*}
Since $2\nu-1+r-p=0$ by assumption, we get 
$v(m,\,m) (\det \bV_{p})^{\nu}\vm=0$. 

Next, let us consider the case that $m < n$. 
By direct computation, we see that 
$v(m,\,n)=(2m)^{-1}[v(m,\,m),\,v(-m,\,n)]_{r}$. 
Since $v(-m,\,n)(\det \bV_{p})^{\nu}\vm=0$ by Case 1, 
and since $v(m,\,m)(\det \bV_{p})^{\nu}\vm=0$ 
by the argument above, it follows that 
$v(m,\,n)(\det \bV_{p})^{\nu}\vm=0$. 
This completes the proof of Proposition~\ref{prop:step4}
\end{proof}
%
%
\subsection{Irreducibility of the quotient modules.}
\label{subsec:quot}
Fix $r \in \BZ$ as above. Denote by $W_{r}^{(1)}$ 
the $\Lon$-submodule of $\Mon$ generated by 
the singular vectors obtained in 
Proposition~\ref{prop:step4}, i.e., 
\begin{equation*}
W_{r}^{(1)}:=\Bigl\langle 
  (\det \bV_{p})^{\nu} \vm \ \Bigm| \ 
  p,\,\nu \in \BZ_{\ge 1}~\text{with}~r=1-2\nu+p
\Bigr\rangle \subset \Mon.
\end{equation*}
%
%
\begin{prn} \label{prop:max3}
A singular vector of $\Mon$ is equal to 
a scalar multiple of the singular vector 
$(\det \bV_{p})^{\nu} \vm$ 
for some $p,\,\nu \in \BZ_{\ge 1}$ with $r=1-2\nu+p$. 
Therefore the $\Lon$-submodule $W_{r}^{(1)}$ 
is the maximal proper submodule of $\Mon$, 
and hence the quotient $\Lon$-module 
$\Mon/W_{r}^{(1)}$ is irreducible. 
\end{prn}

\begin{proof}
Let $u \in \Mon$ be a singular vector, and 
assume that the weight of $u$ is equal to 
$\lambda \in Q_{+}^{(1)} \stm \{0\}$.
By Claim~\ref{c:s3-2} in the proof of 
Proposition~\ref{prop:step3} and 
by Remark~\ref{rem:nu}, we see that the weight $\lambda$ is 
of the form: $2\nu \, \sum_{q=1}^{p}\Lambda_{-q}$ 
for some $p > 0$ and $\nu > 0$, and that 
$w(\bx_{\lambda}) \in \BB[u]$. 
In addition, it follows from 
\eqref{eq:rnup} that $r=2\nu-1-p$. 
Thus, by Proposition~\ref{prop:step4}, 
$(\det \bV_{p})^{\nu} \vm$ is a singular vector.
Because $w(\bx_{\lambda}) \in \BB[(\det \bV_{p})^{\nu} \vm]$ 
by definition, there exists $\alpha \in \BC \stm \{0\}$ 
such that 
%
%
\begin{equation} \label{eq:sg1}
w(\bx_{\lambda}) \notin \BB[u-\alpha(\det \bV_{p})^{\nu} \vm]. 
\end{equation}
Here we should remark that 
$u-\alpha(\det \bV_{p})^{\nu} \vm$ is also 
a singular vector of weight $\lambda$ if it is nonzero. 
Therefore we deduce from \eqref{eq:sg1} and 
Claim~\ref{c:s3-2} in the proof of 
Proposition~\ref{prop:step3} that 
$u-\alpha(\det \bV_{p})^{\nu} \vm=0$, and hence 
$u =\alpha (\det \bV_{p})^{\nu} \vm$. 
Thus we have proved the first assertion of the proposition. 
The other assertions are obvious. 
This completes the proof of the proposition. 
\end{proof}

Let $I_{r}$ be the ideal of the VOA $\Vam \ (=M_{r})$ 
generated by all $(\det \bV_{p})^{\nu} \vm$ for 
$p,\,\nu \in \BZ_{\ge 1}$ with $r=1-2\nu+p$. 

\begin{cor}
The ideal $I_{r}$ is the maximal proper 
ideal of the VOA $\Vam$. Therefore 
the quotient VOA $\Vam/I_{r}$ is simple. 
\end{cor}

\begin{proof}
It follows from Propositions~\ref{prop:max2} 
and \ref{prop:max3} that $W_{r}:=U(\CL_{r})W_{r}^{(1)}$ is 
the maximal proper $\CL_{r}$-module of $\Mam$. 
Then we see by Corollary~\ref{cor:max1} that 
$W_{r}$ is the maximal proper ideal of the VOA $\Vam$.
So, let us show that $I_{r}=W_{r}$. 
The inclusion $I_{r} \subset W_{r}$ follows from 
the fact that the ideal $W_{r}$ contains all 
$(\det \bV_{p})^{\nu} \vm$ 
for $p,\,\nu \in \BZ_{\ge 1}$ with $r=1-2\nu+p$. 
Thus, $I_{r}$ is a proper ideal of $\Vam$, 
which implies that 
$I_{r}$ is a proper $\CL_{r}$-submodule of $\Mam$ 
(see the proof of Proposition~\ref{prop:step1}). 
Since $I_{r}$ contains all $(\det \bV_{p})^{\nu} \vm$ 
for $p,\,\nu \in \BZ_{\ge 1}$ with $r=1-2\nu+p$, 
it can be easily seen from the 
definition of $W_{r}$ that $I_{r} \supset W_{r}$. 
Thus we get $I_{r}=W_{r}$, thereby completing 
the proof of the corollary. 
\end{proof}


{\small
\setlength{\baselineskip}{13pt}
\renewcommand{\refname}{References}

}

\end{document}